\documentclass[preprint,12pt]{elsarticle}

\usepackage{amsfonts}
\usepackage{amsmath}
\usepackage{amssymb}
\usepackage{amsthm}
\usepackage{graphics} %% add this and next lines if pictures should be in esp format
\usepackage{epsfig} %For pictures: screened artwork should be set up with an 85 or 100 line screen
\usepackage{latexsym}
\usepackage{subfigure}

\newtheorem{theorem}{Theorem}
\newtheorem{lemma}{Lemma}

\theoremstyle{definition}

\newcommand{\drop}[1]{}
\newcommand{\no}{\noindent}
\newcommand{\pt}{{\partial_t}}
\newcommand{\eps}{\varepsilon}
\newcommand{\vfi}{\varphi}
\newcommand{\R}{\mathbb{R}}

\newcommand{\fer}[1]{(\ref{#1})}
\newcommand{\qtext}[1]{\quad\text{#1}}
\newcommand{\abs}[1]{| #1 |}
\newcommand{\nor}[1]{\| #1 \|}
\newcommand{\grad}{\nabla}
\def\O{\Omega}

\def\Q{Q_{T}}

\newcommand{\wto}{\rightharpoonup}
\newcommand{\M}{{\mathcal{M}}}
\newcommand{\B}{{\mathcal{B}}}

\newcommand{\esca}[1]{\left\langle #1 \right\rangle}

\DeclareMathOperator{\Div}{div}

\journal{Nonlinear Analysis. Real World Applications}        
        
\begin{document}

\begin{frontmatter}

\title{On a cross-diffusion segregation problem arising from a model of interacting particles 
\tnoteref{t1}}
\tnotetext[t1]{First author supported by the Spanish MEC Project
MTM2010-18427. Second author supported by the Spanish MCINN Project MTM2010-21135-C02-01}

\author{Gonzalo Galiano\corref{cor1}} 
\ead{galiano@uniovi.es}

\author{Virginia Selgas}
\ead{selgasvirginia@uniovi.es}

\address{Dpto.\ de Matem\'aticas, Universidad de Oviedo,
c/ Calvo Sotelo, 33007-Oviedo, Spain}

\cortext[cor1]{Corresponding author. Phone:+34 985103343 Fax: +34 985103354}

\begin{abstract}
We prove the existence of solutions of a cross-diffusion parabolic population problem. The  
system of partial differential equations is deduced as the limit equations satisfied by the densities corresponding to an interacting particles system modeled by stochastic differential equations.
According to the values of the diffusion parameters related to the intra and inter-population repulsion
intensities, the system may be classified in terms of an associated matrix. For proving the existence of solutions when the matrix is  positive definite, we use a fully discrete finite element approximation  
in a general functional setting. If the matrix is only positive semi-definite, we use a regularization technique based on a related
cross-diffusion model under more restrictive functional assumptions. 
We provide some numerical experiments demonstrating the weak and strong segregation effects corresponding
to both types of matrices. 
\end{abstract}

\begin{keyword}
Cross-diffusion system, population dynamics, interacting particles modeling, existence
of solutions, finite element approximation, numerical examples. 
\end{keyword}

% \begin{AMS}
% 35K55, 35D30, 92D25.
% \end{AMS}

\end{frontmatter}

% \pagestyle{myheadings}
% \thispagestyle{plain}
% \markboth{G. GALIANO AND V. SELGAS}{A CROSS-DIFFUSION SEGREGATION PROBLEM}

\section{Introduction}

The effects of spatial cross-diffusion on interacting population models have
been widely 
studied since Kerner \cite{kerner59} and Jorn\'e \cite{jorne77} examined the
linear cross-diffusion 
model 
\begin{equation*}
\pt u_i - a_{i1}\Delta u_1 -a_{i2}\Delta u_2 = (-1)^{i+1} u_i(\alpha_i -\beta_i
u_j), 
\end{equation*}
with non-negative self-diffusivities $a_{ii}$, and non-zero
cross-diffusivities $a_{ij}$, for $i,j=1,2,~i\neq j$, and demonstrated that
while self-diffusion tends to damp out all spatial variations in the
Lotka-Volterra system, cross-diffusion may give rise to instabilities \cite{okubo01}
and to non-constant stationary solutions.

First nonlinear cross-diffusion models seem to have been introduced by 
Busenberg and Travis \cite{busenberg83} (see also Gurtin and Pipkin \cite{gurtin84} for a related model), and
Shigesada et al. \cite{shigesada79} from different modeling points
of view. Shigesada et al. approach starts with the assumption of a single
population density evolution determined by a continuity equation
\begin{equation}
\label{S1}
 \pt u -\Div J(u) =u(\alpha -\beta u), \qtext{with }J(u)=\grad ((c +a u)u) +b u\grad\Phi.
\end{equation}
The divergence of the flow $J$ is thus decomposed into three terms: a random dispersal, $c\Delta u$, 
a dispersal caused by \emph{population pressure}, $a\Delta u^2$, and a drift directed to the
minima of the environmental potential $\Phi$. Generalizing this scalar equation to two 
populations they propose the system, for $i=1,2$,
\begin{equation*}
 \pt u_i -\Div J_i (u_1,u_2)=f_i(u_1,u_2),
\end{equation*}
with
\begin{equation}
\label{S3}
 J_i(u_1,u_2)=\grad \big((c_i +a_{i1}u_1+a_{i2}u_2)u_i\big) +b_i u_i\grad\Phi ,
\end{equation}
and $f_i$ of the competitive Lotka-Volterra type.
Disregarding the linear dispersals ($c=c_i=0$) representing a random contribution to 
the motion, the nonlinear part of the flow $J$ in Eq. \ref{S1} may be 
expressed in conservative form as $J(u)=u\tilde J(u)$, with $\tilde J$ given by the potential 
$\tilde J(u)=\grad (2au+b \Phi)$. However, rewriting the flows \fer{S3} in a similar way leads 
to  the more intricate expression
\begin{equation*}
\tilde J_i(u_1,u_2)= \big(2a_{ii}u_i+a_{ij}\frac{u_j}{u_i}\big)\grad u_i +a_{ij}\grad u_j + b_i \grad\Phi,
\end{equation*}
which, in general, can not be deduced from a potential. This fact has been one of the main
difficulties in finding appropriate conditions ensuring the existence of
solutions to the model proposed by Shigesada et al. (\emph{SKT model}, from now on), see 
\cite{kim84,deuring87,yagi93,lou-ni_96b,ggj01,ggj03,chen04,barret04,wen09,galiano12} and their references.

The generalization of the flow in \fer{S1} to several populations (with $c=b=0$)
given by Busenberg and Travis \cite{busenberg83} is perhaps more natural from the
modeling point of view. They assume that the individual population flow $J_i$
is proportional to the gradient of a potential function, $\Psi$, that only
depends on the total population density $U=u_1+u_2$, 
\begin{equation*}
 J_i(u_1,u_2)=a\frac{u_i}{U}\grad \Psi (U).
\end{equation*}
Note that in this way the flow of $U$ is
still given in the form \fer{S1}, with $J(U)=a\grad\Psi(U)$ (and $c=b=0$).
Assuming the power law $\Psi(s)=s^2/2$, we obtain individual
population flows given by
\begin{equation}
\label{flow_gurtin}
 J_i(u_1,u_2)=a u_i \grad U,
\end{equation}
as those introduced by Gurtin and Pipkin \cite{gurtin84} and mathematically analyzed
by Bertsch et al. \cite{bertsch85,bertsch12}. 

In this article we propose a generalization of the Busenberg-Gurtin model
consisting on the assumption that the individual flows $J_i$ depend, instead of
in the total population density $u_1+u_2$, in a general linear combination of both
population densities, possibly different for each population. As remarked in \cite{gurtin84},
 these weighted sums are motivated when considering a set of species with different 
 characteristics, such as size, behavior with respect to overcrowding, etc. 
 In addition, we also
assume that the flows may contain environmental and random effects, which altogether
lead to the following form 
\begin{equation*}
J_i(u_1,u_2)=u_i  \grad (a_{i1}u_1+a_{i2}u_2 +b_i\Phi) + c_i\grad u_i,
\end{equation*}
 which (for $c_i=0$) has a conservative form similar to that of the scalar case.
 We shall refer to this model as the \emph{BT model}.

Let us finally remark that cross-diffusion parabolic systems have been used 
to model a variety of phenomena ranging from ecology \cite{gilad07,tian10,gv11,galiano12,sherrat00,aly04},
to semiconductor theory \cite{chen07,degond97}, granular materials \cite{aranson02,gjv03,marques03} 
or turbulent transport in plasmas \cite{castillo02},  among others.
Apart from global existence and regularity results for the evolution problem, 
construction of traveling wave solutions \cite{wu05} or exact solutions \cite{cherniha08}
have been accomplished. For the steady state problem, 
existence of non constant steady state solutions has been proven in \cite{lou-ni_96a,lou-ni_96b}. 
Other interesting properties, such as pattern formation, has been 
studied in \cite{vanag09,gambino12,ruiz12,gambino13}. Finally, the numerical discretization has received 
much attention, and several schemes have been proposed \cite{ggj01,ggj03,barret04,gambino09,andreianov11,berres11}.

The article is organized as follows. In Section 2, for a better physical understanding 
of our model, we sketch a heuristic 
deduction based on stochastic dynamics of 
particle systems. In Section 3 we give the
precise 
assumptions on the data problem and state the main results. In Section 4, we introduce the approximated problems
and perform some numerical experiments showing the behavior of solutions under several choices of 
the parameters, including a comparison between the SKT and the BT models. In Section 5, 
we prove the theorems stated in Section 3 , finally, in Section 6 we present our conclusions.

\section{Mathematical modeling}
In recent years there has been a trend to the rigorous deduction of Eq. \fer{S1}
as the equation satisfied by the limit density distribution of suitable particle
stochastic systems of differential equations, see \cite{morale98,
oelschlager89,stevens00,lachowicz11} and their references. We sketch here the
formulation and the main ideas contained in these works which allow us to 
deduce our model.

 Consider a system of $N=N_1+N_2$ interacting particles of two different types 
 described by their trajectories 
$X^i_{j_i}:\R_+\to\R^m $, $j_i=1,\ldots ,N_i$,  $i=1,2$
(stochastic processes). We take $N_1=N_2=n$ to simplify the notation. The
Lagrangian approach to the description of the system is based on
specifying suitable
interacting laws among particles in such a way that their trajectories are
determined by solving 
the following stochastic system of ordinary differential equations (SDE)
\begin{equation}
\label{eq:SDE}
 dX^i_{j}(t)=
F^i_{j}(X^1_1(t),\ldots,X^1_{n}(t),X^2_1(t),\ldots,X^2_{n}(t))dt
+\sigma^i_n dW^i_{j}(t),
\end{equation}
together with some initialization of the processes $X^i_{j}(0)=X^i_{j0}$, 
$j=1,\ldots ,n$, $i=1,2$. 
Functions $F_{j}^i:\R^{2n} \to\R^m $ describe deterministic
interactions among
particles while the constants $\sigma^i_N$
are the intensities of random
dispersal, due to a variety of factors, described by the Brownian motions
$W^i_{j}$,
with $\left\{W^i_1,\ldots,W^i_{n}\right\}$, $i=1,2$,  two families of independent 
standard Wiener
processes valued in $\R^m $.

The individual  particles state may be modeled as positive Radon measures
\begin{equation*}
 \epsilon_{X^i_j(t)}(B)=\left\{
\begin{array}{ll}
 1 & \text{if }X^i_j(t)\in B\\
0 &  \text{if }X^i_j(t)\notin B
\end{array}
\right. \qtext{for all } B\in \B(\R^m ),
\end{equation*}
where $\B(\R^m )$ denotes the Borel $\sigma-$algebra generated by open sets in $\R^m$, while the collective behavior of the discrete system may be given in 
terms of the spatial distribution of particles at time $t$, 
expressed through the \emph{empirical measures}
\begin{equation}
\label{empirical}
 u^i_{n}(t)=\frac{1}{n}\sum_{j=1}^{n} \epsilon_{X^i_{j}(t)}  \in \M(\R^m
),
\end{equation}
which give the spatial relative frequency of particles of the $i$-th population, at time $t$.
Introducing, for 
$\eps>0$, a regularization-scaling kernel $\zeta_\eps(\cdot)=\eps^ {-m}\zeta(\cdot/\eps)$, 
with $\zeta \in C_0^ \infty(\R^ m)$, $\zeta \geq 0$ and $\int\zeta =1$, 
 we may assume that 
the force exerted on the $j$-th single particle of the $i$-th population 
located at $X^i_{j}(t)$ due to the
interaction with all the other particles is given by 
\begin{equation*}
 I^i_j=\sum_{k=1}^2\frac{a_{ik}}{n}\sum_{l=1}^{n}
\zeta_\eps(X^i_j(t)-X^k_l(t)),
\end{equation*}
which may be expressed, using the convolution product, as
\begin{equation*}
I^i_j=\sum_{k=1}^2a_{ik}(u^k_{n}(t)*\zeta_\eps)(X^i_j(t)) .
\end{equation*}
Here, the non-negative coefficients $a_{ik}$ represent the repulsion, pressure or compression intensity
of inter- and intra-specific types, while parameter $\eps$ determines the type of interaction: 
macro, micro or mesoscale, see \cite{morale98}.
The Lagrangian description of the dynamics of our system of interacting
particles
\fer{eq:SDE} may be rewritten in terms of the empirical measures   as
\begin{equation}
\label{eq:SDE2}
 dX^i_{j}(t)= F^i[u^1_{n}(t),u^2_{n}(t)](X^i_{j}(t))dt+
\sigma^i_n dW^i_{j}(t),\quad j=1,\ldots,n.
\end{equation}
We distinguish two kinds of deterministic interactions assuming 
$F^i=F^i_1+F^i_2$, with $F^i_1$ a repulsive interaction between particles given
as
\begin{equation*}
F^i_1[u^1_{n}(t),u^2_{n}(t)](X^i_{j}(t))=-\grad I_j^i=
-\sum_{k=1}^2a_{ik}(u^k_{n}(t)*\grad\zeta_\eps)(X^i_j(t)) , 
\end{equation*}
and $F_2$ a local force, independent of the scaling parameter, 
derived from a potential $\Phi:\R^m \to\R$ 
\begin{equation*}
F_2[u^1_{n}(t),u^2_{n}(t)](X^i_{j}(t))=b_i\grad\Phi (X^i_{j}(t)),
\end{equation*}
with $b_i\in\R$. Finally, with respect to the stochastic part of system
\fer{eq:SDE2}, we assume 
\begin{equation*}
  \lim_{n\to\infty}\sigma^i_n = \sigma_i \geq 0.
\end{equation*}
 Observe that, in some contexts,  $\sigma^i_n$ stands for the \emph{mean
 free path}, i.e., the average distance covered by a moving particle between
 successive collisions. Therefore, the sequence $\sigma^i_n >0$ should be
decreasing with respect to $n$, and a vanishing limit $\sigma_i$ must not be discarded.

\subsection{The Euler description}
A fundamental tool in the derivation of the Eulerian model corresponding to
the Lagrangian description \fer{eq:SDE2} is Ito's formula 
for the time evolution of any smooth scalar function $f(X^i_{j}(t),t)$. 
Introducing the notation
\begin{equation*}
\esca{\mu,g}=\int g(s) \, d\mu(s),
 \end{equation*}
 for the duality $\M(\R^m)\times \mathcal{C}_0(\R^m )$, we deduce  for $i=1,2$
\begin{eqnarray}
 \esca{u^i_{n}(t),f(\cdot,t)}= &&\frac{1}{n}\sum_{j=1}^{n}f(X_j^i(t),t)
=\esca{u^i_{n}(0),f(\cdot,0)} \label{eq:ito_agregada}\\
&& -\sum_{k=1}^2a_{ik} \int_0^t \esca{u^i_{n}(s),
(u^k_{n}(s)*\grad\zeta_\eps)(\cdot) \grad f(\cdot,s )} ds \nonumber\\
&& +b_i \int_0^t \esca{u^i_{n}(s), \grad\Phi \cdot \grad
f(\cdot,s )} ds \nonumber \\
&& + \int_0^t \esca{u^i_{n}( s),\frac{\partial}{\partial s} f(\cdot,s ) +
\frac{1}{2} (\sigma^i_{n})^2  \Delta f(\cdot,s ) }
ds \nonumber\\
&& +\frac{\sigma^i_{n}}{n} \sum_{j=1}^{n} \int_0^t\grad  f(X^i_j(s),s )
\cdot dW^i_j(s) .\nonumber
\end{eqnarray}
The last term of this identity
\begin{equation*}
M^i_{n}(f,t)= \frac{\sigma_{n}^i}{n}\sum_{j=1}^{n} \int_0^t \grad 
f(X^i_j(s),s ) \cdot dW^i_j(s) 
\end{equation*}
is the only explicit source of stochasticity in the equation and shows how,
when the number of particles $n$ is large but still
finite, also from the Eulerian point of view the system keeps the stochasticity
which characterizes each individual. However, Doob’s inequality \cite{friedman75}
 implies that $M^i_{n}(f,t)\to 0$ as $n\to\infty$ in
probability, for any $f\in L^\infty (0,T;W^{1,\infty}(\R^m ))$. In other words,
the Eulerian description becomes deterministic when the size of the
particle system tends to infinity.

Now assume that $u^i_{n}(t)$ tends, as $n\to\infty$, to a deterministic
process
$u^i_\infty(t)$ which 
may be represented by a density function $u_i$ with respect
to the Lebesgue measure on $\R^m $ so that, for  any $t>0$
\begin{equation*}
 \lim_{n\to\infty}\esca{u^i_{n}(t),f(\cdot,t)}=\esca{u^i_\infty(t),f(\cdot,
t) }
=\int_ { \R^m} f(x,t)u_i(x,t)dx.
\end{equation*}
 Then, in the limit $n\to\infty$ we formally obtain from \fer{eq:ito_agregada}
 (see \cite{capasso08,lachowicz11} for the rigorous deduction of this limit)
\begin{eqnarray*}
\int_{\R^m } f(x,t)u_i(x,t)dx=&& \int_{\R^m } f(x,0)u_i(x,0)dx \\
&&-\sum_{k=1}^2a_{ik} \int_0^t\int_{\R^m } u_i(x,s) \grad u_k(x,s)\cdot \grad
f(x,s) dx ds\\
&&+b_i \int_0^t\int_{\R^m } u_i(x,s) \grad\Phi(x)\cdot \grad f(x,s) dx ds\\
&&+ \int_0^t\int_{\R^m } u_i(x,s)\big( 
\frac{\partial}{\partial s} f(x,s) +
\frac{1}{2} \sigma_i^2  \Delta f(x,s ) \big) dx ds,
\end{eqnarray*} 
which may be recognized as a weak formulation of the following Cauchy PDE
problem for the unknowns $u_i:\R^m \times\R_+\to [0,1]$
\begin{equation}
\label{cauchy}
\pt u_i-\Div\big(u_i(a_{i1}\grad u_1 +a_{i2}\grad u_2 -b_i\grad \Phi)\big) -
c_i\Delta u_i=0
\qtext{in }\R^m\times\R_+,  
\end{equation}
for initial data $u_i(\cdot,0)=u_{i0}$ in $\R^m $, and $c_i=\sigma_i^2/2$. 

Let us, finally, remark that the deduction of the Cauchy problem \fer{cauchy} is not easily
extended to boundary value problems. In which respects to the non-flow boundary conditions studied in the
next section, the corresponding SDE system seem to be the so-called Skorohod or reflecting boundary (stochastic) problem, in which
particles are reflected in some prescribed direction when hitting the boundary. Although there exists an abundant literature on this 
problem, see for instance \cite{lions84,slominski93,singer08} and their references, to the knowledge of the authors there is not a rigorous 
deduction of a PDE problem satisfied by the corresponding limit density.

\section{Assumptions and main results}
Inspired by the problem deduced in the previous section, we set the following one: Given a fixed $T>0$ and a bounded set $\O\subset\R^m$,  find
$u_i:\O\times(0,T)\to\R$ such that, for $i=1,2$,
\begin{align}
& \pt u_i-\Div J_i(u_1,u_2)=f_i(u_1,u_2) && \qtext{in }Q_T=\O\times(0,T), 
  \label{eq:pde}\\
& J_i(u_1,u_2)\cdot n =0 && \qtext{on }\Gamma_T=\partial\O\times(0,T),
\label{eq:bc}\\
& u_{i}(\cdot,0)=u_{i0} && \qtext{in }\O,	\label{eq:id} 
\end{align}
with flow and competitive Lotka-Volterra functions given by
\begin{align}
 & J_i(u_1,u_2) = u_i\big(a_{i1}\grad u_1 +a_{i2}\grad u_2 +b_i q \big) +
 c_i\grad u_i ,  \label{def:flow}\\
& f_i(u_1,u_2) = u_i \big(\alpha_{i}-\beta_{i1} u_1 - \beta_{i2} u_2 \big) ,
\label{def:reaction}
\end{align}
where the coefficients $a_{ij},~c_i,~b_i,~\alpha_i,~\beta_{ij}$, $i,j=1,2$ are assumed to be functions, and not merely constants. Observe that we also replaced the potential field $\grad \Phi$ of the model derived in the previous section by a general field $q$.
We make the following hypothesis on the data, which we shall refer to as \textbf{(H)}:
\begin{enumerate}
\item $\O\subset\R^m$ ($m=1,2$ or $3$) is a bounded set with Lipschitz continuous boundary $\partial \O$.
\item For $i,j=1,2$,  the coefficients $a_{ij},~c_i,~\alpha_i,~\beta_{ij} \in L^{\infty}(Q_T)$ are non-negative a.e. in $Q_T$, and $b_i\in L^\infty(Q_T)$.
Besides, there exists a constant $a_0 >0$ such that
\begin{equation}\label{H:def_pos}
4 a_{11} a_{22} - (a_{12}+a_{21})^2 > a_0 \qtext{a.e. in } Q_T . 
\end{equation}
\item The drift function satisfies $q \in (L^2(Q_T))^m$.
\item The initial data are non-negative and satisfy $u_{i0}\in L^\infty(\O)$,
$i=1,2$.
\end{enumerate}

Notice that condition \fer{H:def_pos} implies the following ellipticity condition on the matrix $A=(a_{ij})_{i,j=1}^2$:
\begin{equation}\label{cond:ellipticity}
\xi^T A~\xi \geq a_0 \nor{\xi}^2 \qtext{a.e. in }Q_T \text{ and for all } \xi\in\R^2.
\end{equation}

\begin{theorem}\label{th:existence_original}
Let $T>0$ and assume (H). Then problem \fer{eq:pde}-\fer{eq:id} has a weak solution $(u_1,u_2)$ satisfying $u_i\geq 0$ in $Q_T$ and 
$$
u_i\in L^2(0,T;H^1(\O))\cap L^{r} (Q_T) \cap W^{1,p}(0,T;(W^{1,p'}(\O))'),\quad
i=1,2,
$$
where $p=(2m+2)/(2m+1)$, $r=2(m+1)/m$ and $p'=2(m+1)$, in the sense that for all $\vfi\in L^{p'}(0,T;W^{1,p'}(\O))$, $i=1,2$,
\begin{align}\label{weak}
\int_0^T<\pt u_i,\vfi>  
%& +\int_{Q_T} \big(u_i (a_{i1}\grad u_1 +a_{i2}\grad u_ 2 + b_i\bq )+ c_i \grad u_i\big)\cdot\grad\vfi  \\
+\int_{Q_T} J_i(u_1,u_2) \cdot\grad\vfi  %\\ &
= \int_{Q_T} f_i(u_1,u_2)\,\vfi, 
\end{align}
with $<\cdot,\cdot>$ denoting the duality product between 
$W^{1,p'}(\O)$ and its dual $(W^{1,p'}(\O))'$.
\end{theorem}
As in \cite{ggj03} (for $m=1$) and in Chen and J\"ungel \cite{chen04} (for $m\leq 3$), the main tool for the analysis of problem \fer{eq:pde}-\fer{eq:id} is the use of the entropy functional 
\begin{equation}
\label{def:entropy}
E(t)= \sum_{i=1}^2\int _{\Omega} F(u_i(\cdot,t)) \geq 0, \qtext{with }F(s)=s (\ln s-1)+1,
\end{equation}
which allows us to deduce formally the identity 
\begin{align*}
E(t)+ \int_{Q_t } \Big( \sum_{i=1}^2 ( a_{ii} |\grad u_i|^2 + 2c_i\abs{\grad
\sqrt{u_i}}^2 )+ (a_{12}+a_{21}) \grad u_1\cdot \grad u_2 \Big) \\
  = E(0)+
 \int_{Q_t } \sum_{i=1}^2  \Big( -b_i \, q \cdot \grad u_i +
 f_i(u_1,u_2) \ln u_i \Big) , %\qtext{in }(0,T) .
\end{align*}
by using $F'(u_i)=\ln u_i$ as a test function in the weak formulation of \fer{eq:pde}-\fer{eq:id}. 
From assumptions (H) and, specially, bound \fer{H:def_pos}, one easily obtains the entropy inequality
\begin{align}\label{ineq:entropy}
E(t)+ a_0 \int _{Q_t} ( |\grad u_{ 1}|^2 + |\grad u_{ 2}|^2)\leq (E(0)+C_1) \, \mathrm{e} ^{C_2 t} ,
\end{align}
providing the key  $L^2(0, T;H^1(\O­))$ estimate of $u_1$ and $u_2$ which allows to prove 
Theorem~\ref{th:existence_original}.  
Thus,  bound \fer{H:def_pos} provides a sufficient condition on the diffusion operator to prove the
existence of solutions of problem \fer{eq:pde}-\fer{eq:id} under conditions (H). However, these conditions are not necessary, 
as the following result shows. First, we state the precise assumptions to treat this degenerate case, 
to which we refer to as \textbf{(H')}:
\begin{enumerate}
\item The boundary $\partial \O$ is $ H_2$ (H\"older continuous with exponent $2$). 

\item $a_{ij}=a,~b_i=b,~c_i=c,~f_i(s_1,s_2 )= s_i (\alpha -\beta(s_1+s_2))$ for some constants
$a,~b,~c,~\alpha,~\beta$ such that $a >0$ and $c,~\alpha,~\beta \geq 0$.

\item The drift function satisfies $q \in (L^\infty(Q_T))^m$, $\Div q \in L^\infty(Q_T)$.

\item $u_0=u_{10}+u_{20} \in H^2(\O)$,  with  $u_{0} > \tilde u$, for some constant $\tilde u>0$,
and $\big(a u_0 \grad u_0 +b q u_0 +c\grad u_0\big)\cdot n =0$ on $\partial\O$ (compatibility condition).

\end{enumerate}
Under assumptions (H'), the equation satisfied by $u_i$, for $i=1,2$, is 
\begin{equation}
 \pt u_i -\Div \big(a u_i \grad (u_1+u_2) +b q u_i +c\grad u_i\big)= u_i (\alpha -\beta (u_1+u_2)), 
\end{equation}
which is closely related to the model introduced by Gurtin and Pipkin \cite{gurtin84}. An important case included in assumptions (H') is the \emph{contact inhibition problem} arising in tumor modeling, see for instance Chaplain et al. \cite{chaplain06}, i.e. that in which the initial data, describing the spatial distribution of normal and 
tumor tissue,  satisfy 
$\{u_{10}>0\}\cap \{u_{20}>0\} =\emptyset$. This free boundary problem was 
mathematically analyzed by Bertsch et al. 
for one \cite{bertsch85} and several spatial dimensions \cite{bertsch12} by using regular Lagrangian flow techniques. 
However, our approach is different and more general in some aspects, like that of the data regularity assumptions or the consideration of a drift term.
\begin{theorem}\label{th:existence_particular}
Let $T>0$ and assume (H').
Then problem \fer{eq:pde}-\fer{eq:id} has a weak solution $(u_1,u_2)$ such that
\[
 u_i\in L^{\infty} (Q_T) \cap H^{1}(0,T;(H^{1}(\O))'),\quad u_1+u_2 \in L^2(0,T;H^{2}(\O))
\]
and an identity similar to \fer{weak} is satisfied for all $\vfi\in L^2(0,T;H^1(\O))$. 
\end{theorem}
We finish this section by showing some connections between the Shigesada et al. model (SKT) and the Busenberg and Travis model (BT) studied in this article. Let 
\begin{equation}
\label{def.NLflows}
J_i^{BT}(u_1,u_2) = u_i(a_{i1}\grad u_1 +a_{i2}\grad u_2 ), \quad  
 J_i^{SKT}(u_1,u_2) = \grad \big(u_i(a_{i1} u_1 +a_{i2} u_2) \big),
\end{equation}
be the nonlinear diffusive flows corresponding to the BT \fer{def:flow}, and SKT \fer{S3}
models, respectively. First, we observe that 
\begin{equation*}
 J_i^{SKT}(u_1,u_2) =J_i^{BT}(u_1,u_2) + (a_{i1} u_1 +a_{i2} u_2)\grad u_i,
\end{equation*}
indicating that the support of diffusion for $J_i^{SKT}$ is, at least, equal to that of $J_i^{BT}$, and 
explaining the smoother behavior of solutions corresponding to $J_i^{SKT}$ observed in the numerical experiments.
We may approximate $J_i^{BT}$ by introducing the perturbation
\begin{equation}
\label{def.NLflowapp}
 J_i^{BT,\delta}(u_1,u_2) =J_i^{BT}(u_1,u_2) +\frac{\delta}{2}J_i^{SKT}(u_1,u_2),
\end{equation}
for some $\delta>0$. Although $J_i^{BT,\delta}$ can not be recast in the same functional form as $J_i^{SKT}$, 
the diffusion matrices corresponding to both flows share an important property, e.g. both give rise to 
a positive definite matrix once the change of unknowns $u_i=\exp(w_i)$ is introduced. Being this idea the main 
ingredient introduced in \cite{ggj03} for the proof of existence of solutions of the SKT model, 
we may follow the steps given in Chen and J\"ungel \cite{chen04} for proving the existence of 
solutions $(u_1^{(\delta)},u_2^{(\delta)})$ corresponding to the problem with nonlinear flow  
$J_i^{BT,\delta}$
and, after obtaining suitable a priori estimates, pass to the limit $\delta\to 0$ to 
deduce the existence of solutions of problem  \fer{eq:pde}-\fer{eq:id} according to 
conditions (H) or (H'). Although we have followed this approach for the proof of Theorem~\ref{th:existence_particular},
we have preferred to use a direct technique  to prove Theorem~\ref{th:existence_original} 
by adapting the Finite Element Method employed by Barrett and Blowey \cite{barret04} which provides 
a convergent fully discrete numerical scheme for our numerical experiments.

% \hspace{3cm}
% 
% Our proof of Theorem~\ref{th:existence_particular} is based on exploiting the properties of the  
% solution to the problem satisfied  by $u=u_1+u_2$,
% \begin{align}
% & \pt u-\Div \big(a u \grad u +b q u +c\grad u\big)= u (\alpha -\beta u) && \qtext{in }Q_T=\O\times(0,T), 
%   \label{eq:pde_p}\\
% & \big(a u \grad u +b q u +c\grad u\big)\cdot n =0 && \qtext{on }\Gamma_T=\partial\O\times(0,T),
% \label{eq:bc_p}\\
% & u(\cdot,0)=u_{0} && \qtext{in }\O,	\label{eq:id_p} 
% \end{align}
% which is a uniformly parabolic problem due to (H')$_2$ and (H')$_4$. Due to these additional assumptions it may be seen that 
% a suitable $\eps-$regularization of problem \fer{eq:pde}-\fer{eq:id} posses a solution for which weak convergence
% for $u_i^{(\eps)}$, $i=1,2$, as well as strong convergence for $u_1^{(\eps)}+u_2^{(\eps)}$ may be shown in appropriate functional spaces such that the
% limit of the solutions of the regularized problem converges to a solution of problem \fer{eq:pde}-\fer{eq:id}.
% The regularized problems are of a non-trivial nature and strongly related to the Shigesada et al. model for which
% existence of solutions was proven in \cite{ggj03} for one spatial dimension and in Chen and J\"ungel \cite{chen04} for up to three
% spatial dimensions.

\section{Approximated problems and numerical experiments}

In this section we describe the regularized problems and the discretization employed to perform the numerical experiments. For approximating problem \fer{eq:pde}-\fer{eq:id} under conditions (H) we adapted the FEM 
technique used in \cite{barret04}. This FEM approach is also used to discretize the SKT model, i.e.  problem \fer{eq:pde}-\fer{eq:id} with $J_i^{BT}$ replaced by $J_i^{SKT}$, see \fer{def.NLflows}, for comparison purposes.
In Experiments 1 and 2, we show these comparisons for data problem taken from \cite{ggj01}. In general terms, 
the qualitative behavior of solutions is similar, although we may observe that model BT produces less regular solutions than model SKT. Although we lack of a rigorous proof, it seems that solutions of the BT model generate spatial niches. 

For approximating problem \fer{eq:pde}-\fer{eq:id} under conditions (H') we proceed as mentioned in the previous section.
We first replace $J_i^{BT}$ by $J_i^{BT,\delta}$, see \fer{def.NLflowapp}, which has similar structural properties than 
the flow of the SKT model. Then, we use the same approach as that under conditions (H), and inspect the behavior of solutions when
$\delta\to 0$. In Experiments 3 and 4 we present results related to the contact inhibition problem. The most interesting phenomenon is the development of discontinuities of $(u_1^{(\delta)},u_2^{(\delta)})$ in
the contact point as  $\delta\to 0$, indicating a parabolic-hyperbolic transition in the behavior of solutions, as already noticed in \cite{bertsch85}. 

Since the numerical scheme is common for the three nonlinear diffusion flows under study, we describe it 
for the general flow
\begin{equation}
 J^G_i (u_1,u_2,\grad u_1,\grad u_2) = J^D_i(u_1,u_2)+b_iu_iq  + c_i \grad u_i,
\end{equation}
with $D=BT,~SKT$ or $D=BT,\delta$. For the numerical experiments, we have chosen constant coefficients for both
the flow and the Lotka-Volterra terms  and an affine 
environmental field $q$. However, general $L^\infty(Q_T)$ coefficients and $L^2(Q_T)$ environmental field 
may be also considered.
For the time discretization, we take in the experiments a uniform partition of $[0,T]$ of time step $\tau$.
For $t=t_0=0$, set $u_{\eps i}^0=u_i^0$.
Then, for $n\geq 1$ find $u_{\eps i}^{n}$ such that for 
$i=1,2$,  
\begin{equation}\label{eq:pde_discr.s4}
\begin{array}{l}
\frac{1}{\tau}\big( u^n_{\eps i}-u^{n-1}_{\eps i} , \chi )^h
+ \big( J^G_i(\Lambda _{\eps } (u^n_{\eps 1}),\Lambda _{\eps } (u^n_{\eps 2}),\grad u^n_{\eps 1},\grad u^n_{\eps 2}  ) ,\grad\chi \big)^h =\\ [2ex]
\hspace*{1cm} = \big(\alpha_{i} u^n_{\eps i} - \lambda _{\eps } (u^n_{\eps i})
( \beta_{i1} \lambda _{\eps } (u^{n-1}_{\eps 1}) + \beta_{i2} \lambda
_{\eps } (u^{n-1}_{\eps 2}) ) , \chi \big)^h , %\qquad\forall \chi\in S^h .
\end{array}
\end{equation}
for every $ \chi\in S^h $, the finite element space of piecewise $\mathbb{P}_1$-elements. Here, $(\cdot,\cdot)^h$ 
stands for a discrete semi-inner product on $\mathcal{C}(\overline{\Omega} )$.
The parameter $\eps>0$ makes reference 
to the regularization introduced by functions $\lambda_\eps$ and $\Lambda_\eps$, which converge to the identity 
as $\eps\to0$. See the Appendix for details

Since \fer{eq:pde_discr.s4} is a nonlinear algebraic problem, we use a fixed point argument to 
approximate its solution,  $(u_{\eps 1}^n,u_{\eps 2}^n)$, at each time slice $t=t_n$, from the previous
approximation $u_{\eps i}^{n-1}$.  Let $u_{\eps i}^{n,0}=u_{\eps i}^{n-1}$. 
Then, for $k\geq 1$ the problem is to find $u_{\eps i}^{n,k}$ such that for 
$i=1,2$, and for all $\chi \in S^h$ 
\begin{equation*}%\label{eq:pde_discr}
\begin{array}{l}
 \frac{1}{\tau}\big( u^{n,k}_{\eps i}-u^{n-1}_{\eps i} , \chi )^h
+ \big( J^G_i(\Lambda _{\eps } (u^{n,k-1}_{\eps 1}),\Lambda _{\eps } (u^{n,k-1}_{\eps 2}),\grad u^{n,k}_{\eps 1},\grad u^{n,k}_{\eps 2}  ) ,\grad\chi \big)^h =\\ [2ex]
\hspace*{1cm} = \big(\alpha_{i} u^{n,k}_{\eps i} - \lambda _{\eps } (u^{n,k-1}_{\eps i})
( \beta_{i1} \lambda _{\eps } (u^{n-1}_{\eps 1}) + \beta_{i2} \lambda
_{\eps } (u^{n-1}_{\eps 2}) ) , \chi \big)^h . %\qquad\forall \chi\in S^h .
\end{array}
\end{equation*}
We use the stopping criteria $\max _{i=1,2} \nor{u_{\eps,i}^{n,k}-u_{\eps,i}^{n,k-1}}_\infty <\text{tol}$,
% \begin{equation*}
%  \max _{i=1,2} \nor{u_{\eps,i}^{n,k}-u_{\eps,i}^{n,k-1}}_\infty <\text{tol},
% \end{equation*}
for values of $\text{tol}$ chosen empirically, and set $u_i^n=u_i^{n,k}$.
In some of the experiments we integrate in time until a numerical stationary solution, 
$u_i^S$, is achieved. This is determined 
by $\max _{i=1,2} \nor{u_{\eps,i}^{n,1}-u_{\eps,i}^{n,0}}_\infty <\text{tol}_S$,
% \begin{equation*}
%  \max _{i=1,2} \nor{u_{\eps,i}^{n,1}-u_{\eps,i}^{n,0}}_\infty <\text{tol}_S, 
% \end{equation*}
where $\text{tol}_S$ is chosen also empirically. 
Finally, for the spatial discretization we take  a uniform partition of the interval $\O$ in $M$ subintervals.

% \begin{table}
% \centering
% \caption{Parameter values for experiments} 
% \begin{tabular}{|l|c|c|c|}
% \hline\hline 
% Parameter & Symbol & Experiment 1 & Experiment 2  \\ \hline\hline
% Spatial domain & $\O$ & $(0,3)$ & (0,1)\\
% Diffusion coeff. & $(c_1,c_2)$ & $(1,1)$ &  $(0,0)$ \\
% %Cross-diffusion coeff. & $a_{ij}$ & $\left(\begin{array}{cc} 1& 1 \\ 1 & 1\end{array}\right)$ &
% %$\left(\begin{array}{cc} \cdot & 1 \\ 1 & \cdot \end{array}\right)$  \\
% Drift coeff. & $(b_1,b_2)$ & $(\cdot,1)$ & $(0,0)$ \\
% Initial densities & $(u_{10},u_{20})$ & $(10,10)$ & $ (g(x;0.4), g(x;0.6)) ^{(*)}$\\
% Environmental potential & %$\bq$ 
% $q$ & $-3(x-0.5)$ & $0$ \\
% \hline\hline
% \end{tabular}
% \label{table.parametros}
% %\caption{Some typical parameter values} 
% \\
% {\small (*)  $g(x;x_0)=\exp(-\frac{(x-x_0)^2}{0.001})$}
% \end{table}

%\bigskip 

\begin{figure}[t]
\centering
 \subfigure%
 %{\includegraphics[width=6.25cm,height=4.75cm]{exp0a_fig2}}
 {\includegraphics[width=6.25cm,height=4.75cm]{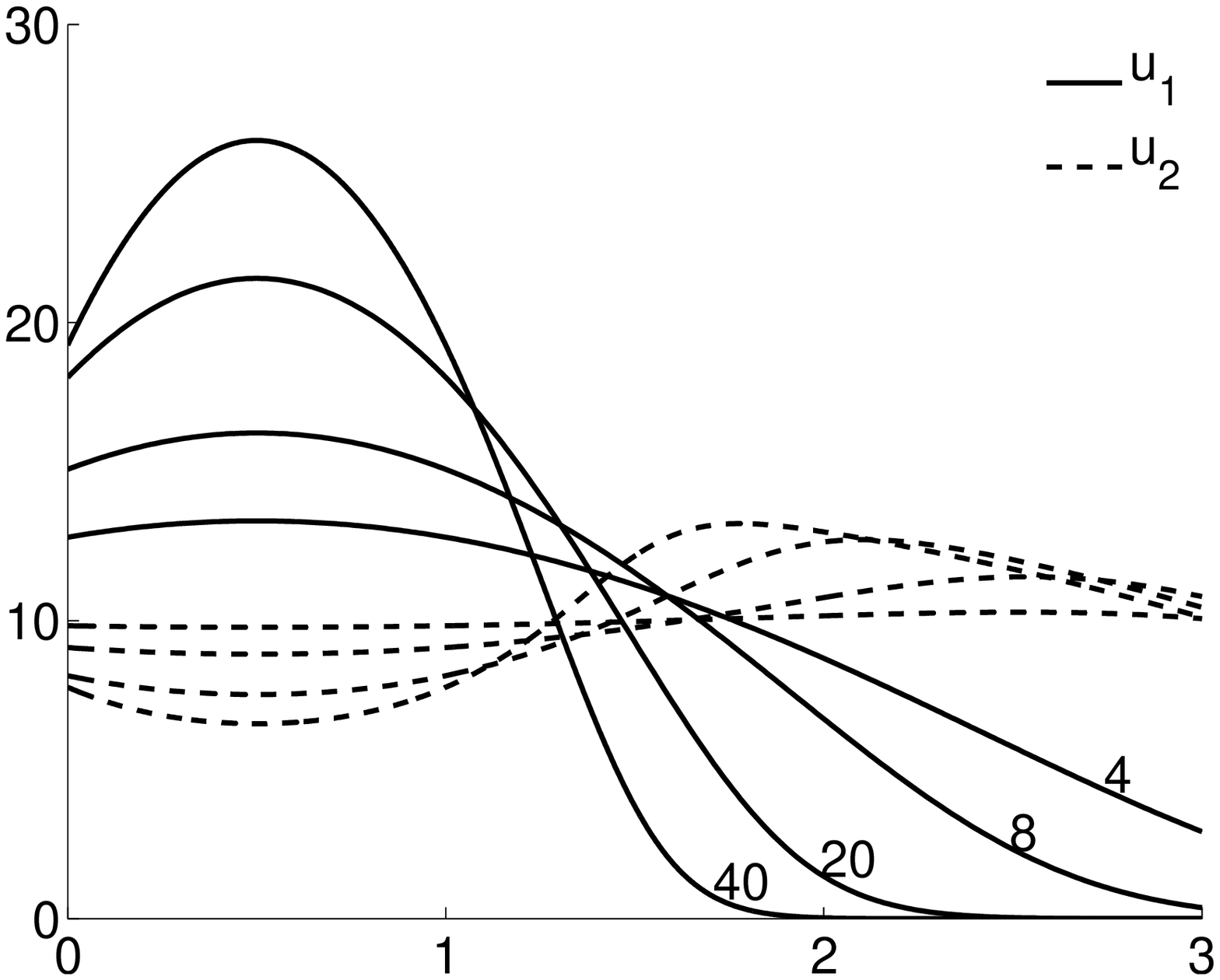}}
 %\hspace{0.5cm}
 \subfigure
 %{\includegraphics[width=6.25cm,height=4.75cm]{exp1a_fig2}} %\\
 {\includegraphics[width=6.25cm,height=4.75cm]{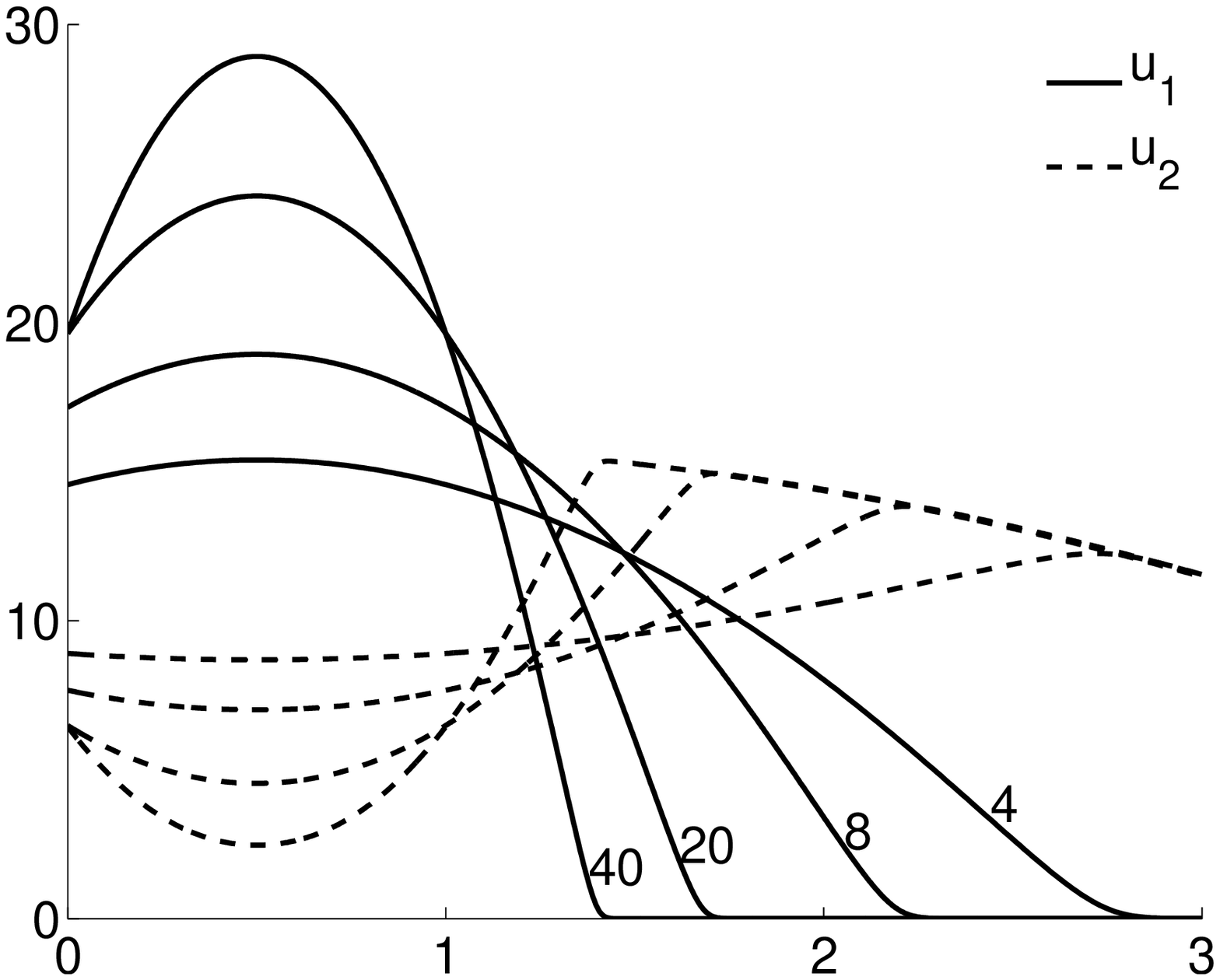}} %\\
%\vspace*{-1cm}
 \caption{{\small Experiment 1. Numerical results for several values of $b_1$ using Shigesada et al. model (left panel)  and the model studied in this article (right panel).}} 
\label{exp01_fig}
\end{figure}

\bigskip

\no\emph{Experiment 1.} We compare graphically the phenomenon of segregation of populations arising
from the Shigesada et al. model \fer{S3}, and from the model studied in this article \fer{eq:pde}-\fer{eq:id}. To do this we use Example 
(c) of \cite{ggj01} in which an implicit finite differences method was used to compute the approximated solution,
see also \cite{barret04,gambino09} for the same experiment reproduced with alternative methods. 
The parameters are fixed as follows: $\O=(0,3)$,  $a_{ij}=1$, $c_i=1$, $i,j=1,2$, $b_2=1$ and $b_1=4,~8,~20,~40$. 
The initial data is constant, $u_{i0}=10$, for $i=1,2$, and the environmental field is given by $q(x)=-3(x-0.5)$.
For the Shigesada et al. model we take $M=301$ and $\tau=10^ {-3}$, as in \cite{ggj01}. However, the convergence properties of problem \fer{eq:pde}-\fer{eq:id} lead us  to choose values of $M$ and $\tau$ in the range $301-1001$ and $10^ {-5}-10^ {-3}$, 
respectively, depending on the $b_1$ values. For both models we take tolerances $\text{tol}\sim 10^{-7}$ and $\text{tol}_S\leq 5\times 10^{-8}$. 

\begin{figure}[t]
\centering
 %\subfigure[$t=0$ ]
 {\includegraphics[width=4.25cm,height=2.7cm]{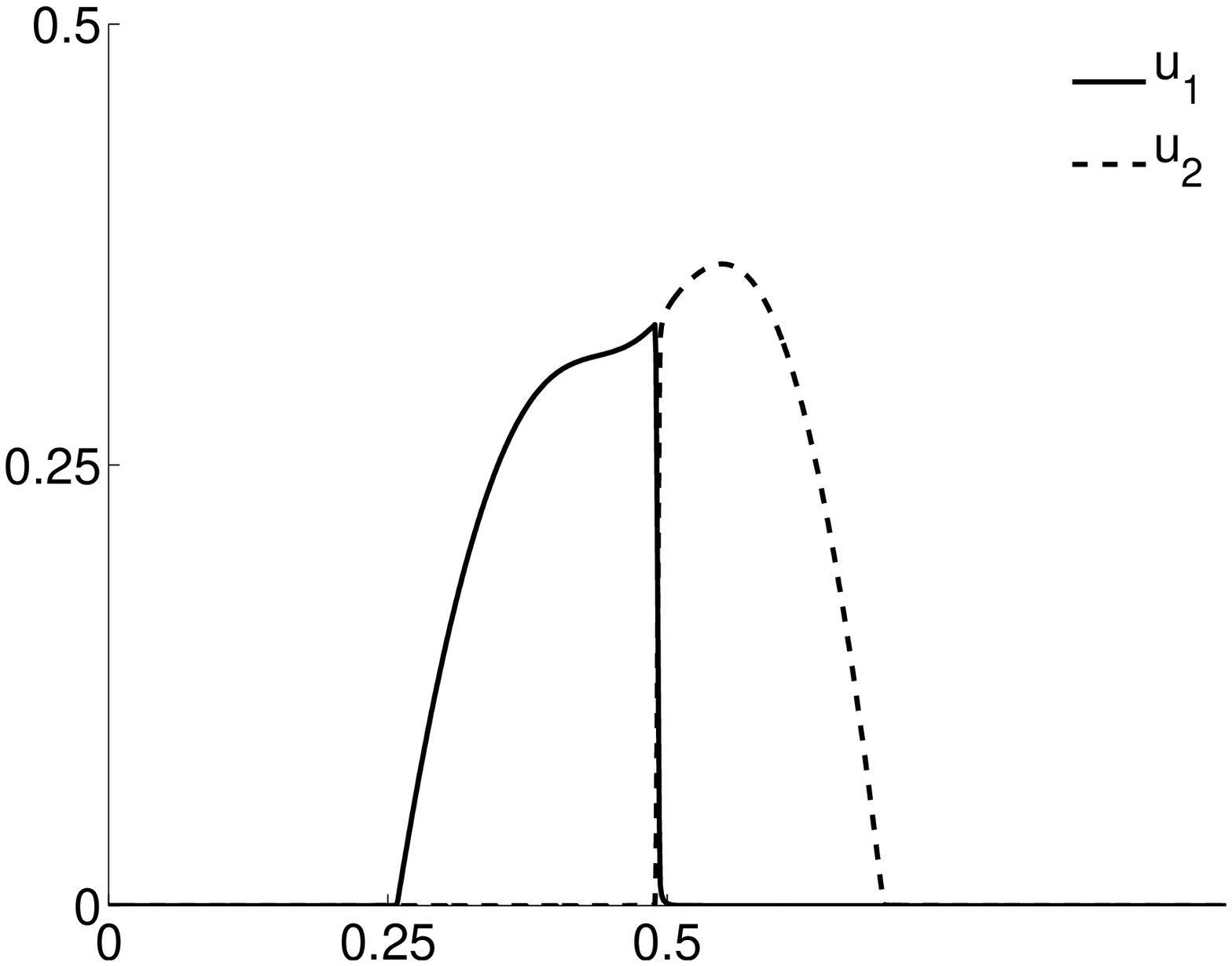}}
 %\hspace{0.5cm}
 %\subfigure[$t=5$]
 {\includegraphics[width=4.25cm,height=2.7cm]{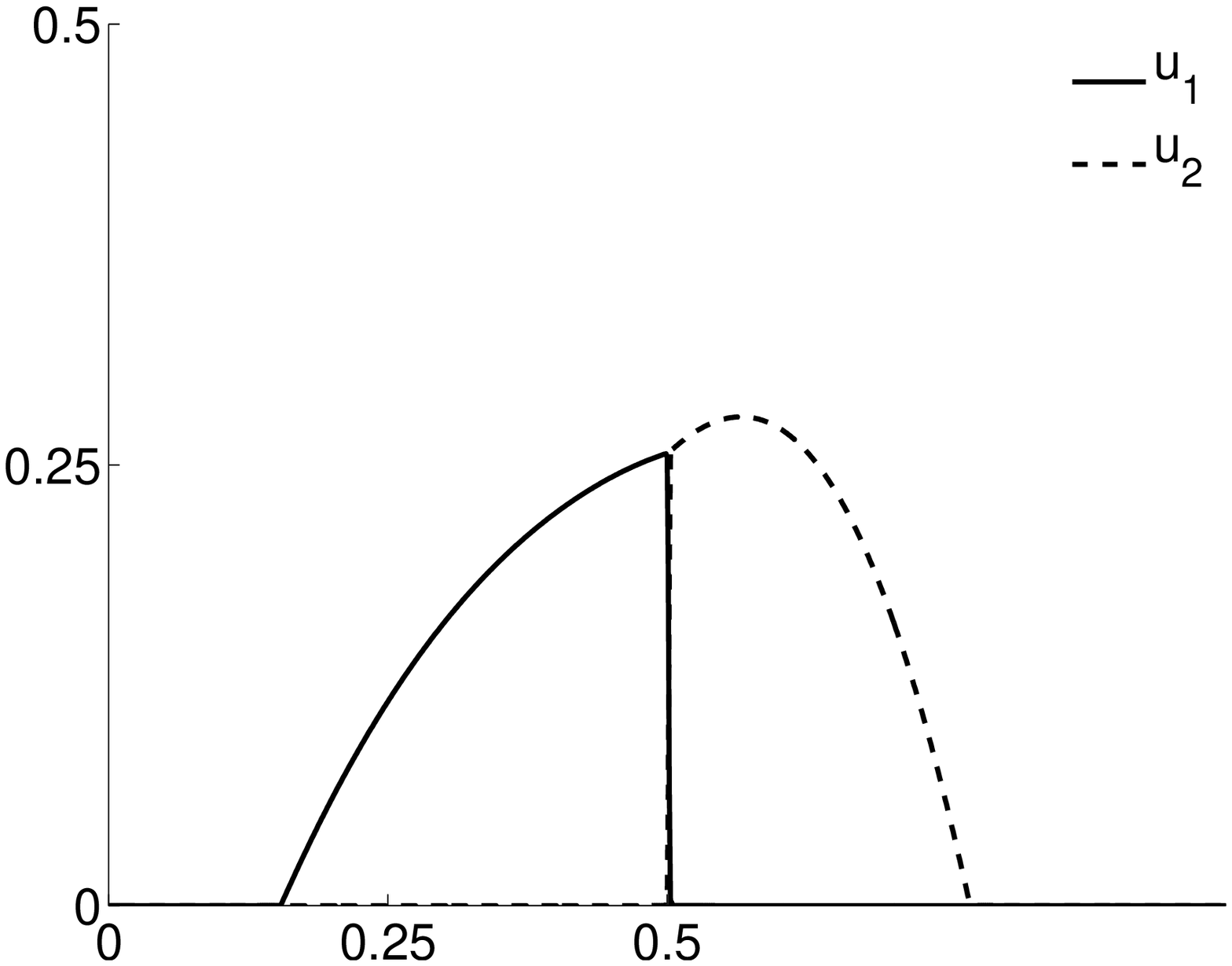}} 
 %\subfigure[$t=7$]
 {\includegraphics[width=4.25cm,height=2.7cm]{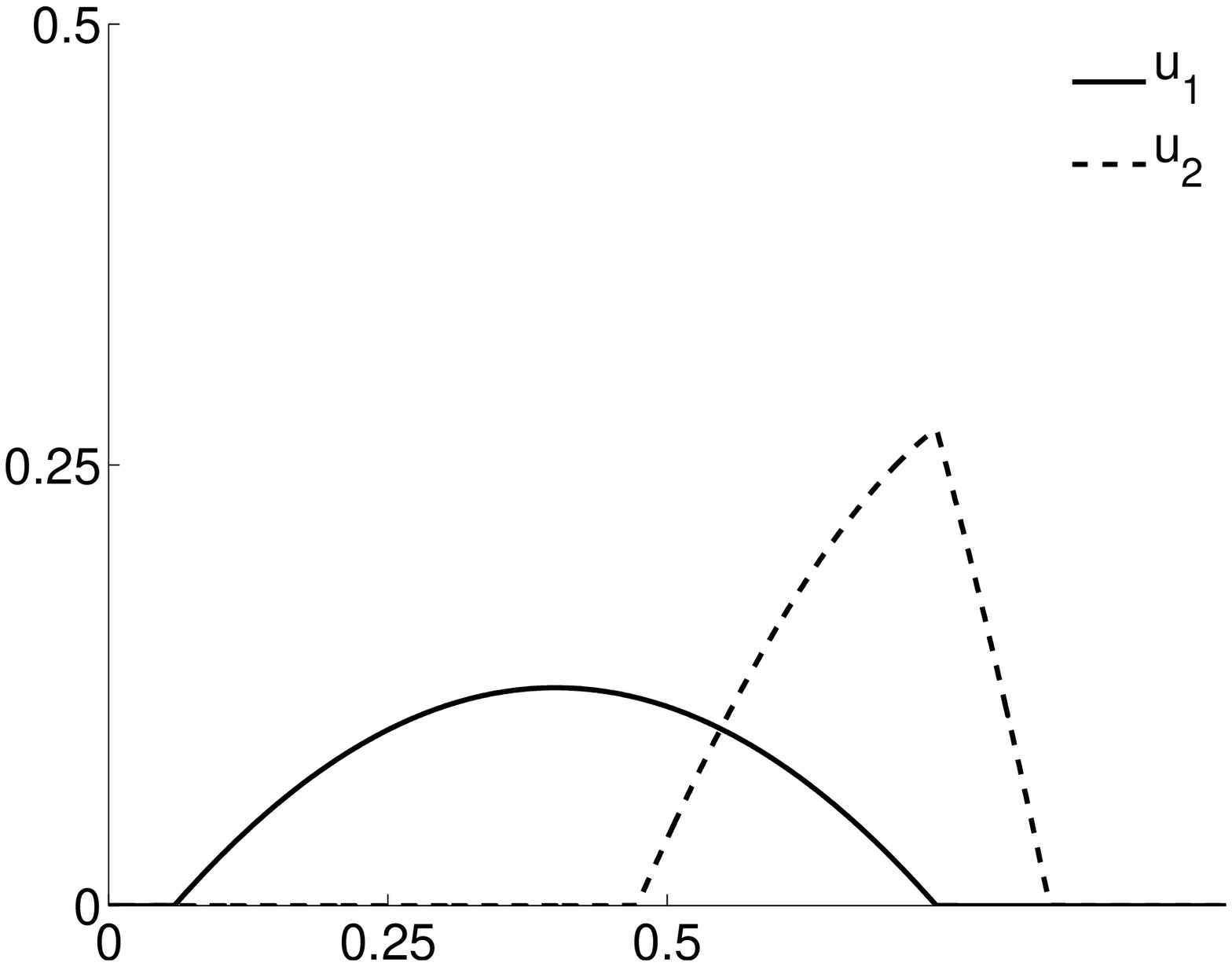}}
%  \hspace{0.5cm}
%  \subfigure[$\delta=0.1$]
%  {\includegraphics[width=6.25cm,height=4.75cm]{exp3e_fig}}
   \caption{{\small Experiments 4-6. Transient state profiles of solutions corresponding to:
   left: $d_1=1,~d_2=10$, center: $a_{11}=a_{12}=3$, $a_{21}=a_{22}=1$, right:
   $a_{11}=4,~a_{12}=0,~a_{21}=3.9,~a_{22}=1$.}} 
\label{exp789_fig}
\end{figure}

In Fig.~\ref{exp01_fig} we plot the approximate steady state solution for both models. Labels on curves give the corresponding $b_1$ value. We observe a stronger segregation effect in the model studied in this article compared to the model of Shigesada et al., although they behave similarly from a qualitative point of view. The loss of regularity 
of solutions when one of them vanishes is also observed. To check this fact more clearly, we run an experiment 
for the same data as above but: $a_{11}=4,~a_{12}=0,~a_{21}=3.9,~a_{22}=1$, 
$b_i=d_i=0$, for $i=1,2$. Observe that matrix $(a_{ij})$ satisfies the positiveness condition 
\fer{H:def_pos}. A transient state of the solution is shown in the right panel of Fig.~\ref{exp789_fig}.

\begin{figure}[t]
\centering
 \subfigure[$t=0$]%
 {\includegraphics[width=4.25cm,height=2.7cm]{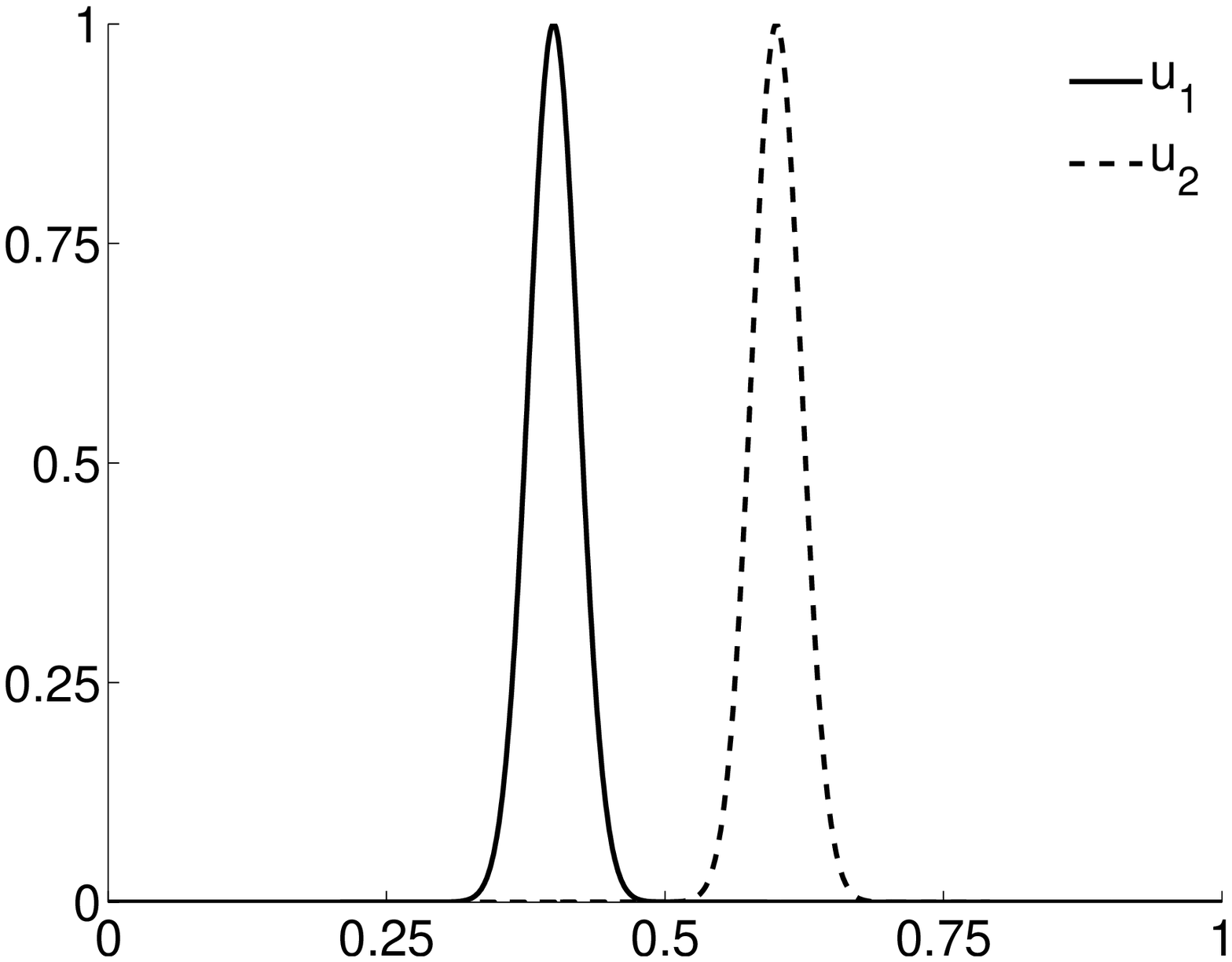}}
  %\hspace{0.5cm}
%  \subfigure[$t=0.01$]
%  {\includegraphics[width=6.25cm,height=4.5cm]{exp2b_fig}} \\
  \subfigure[$t=0.05$]%
 {\includegraphics[width=4.25cm,height=2.7cm]{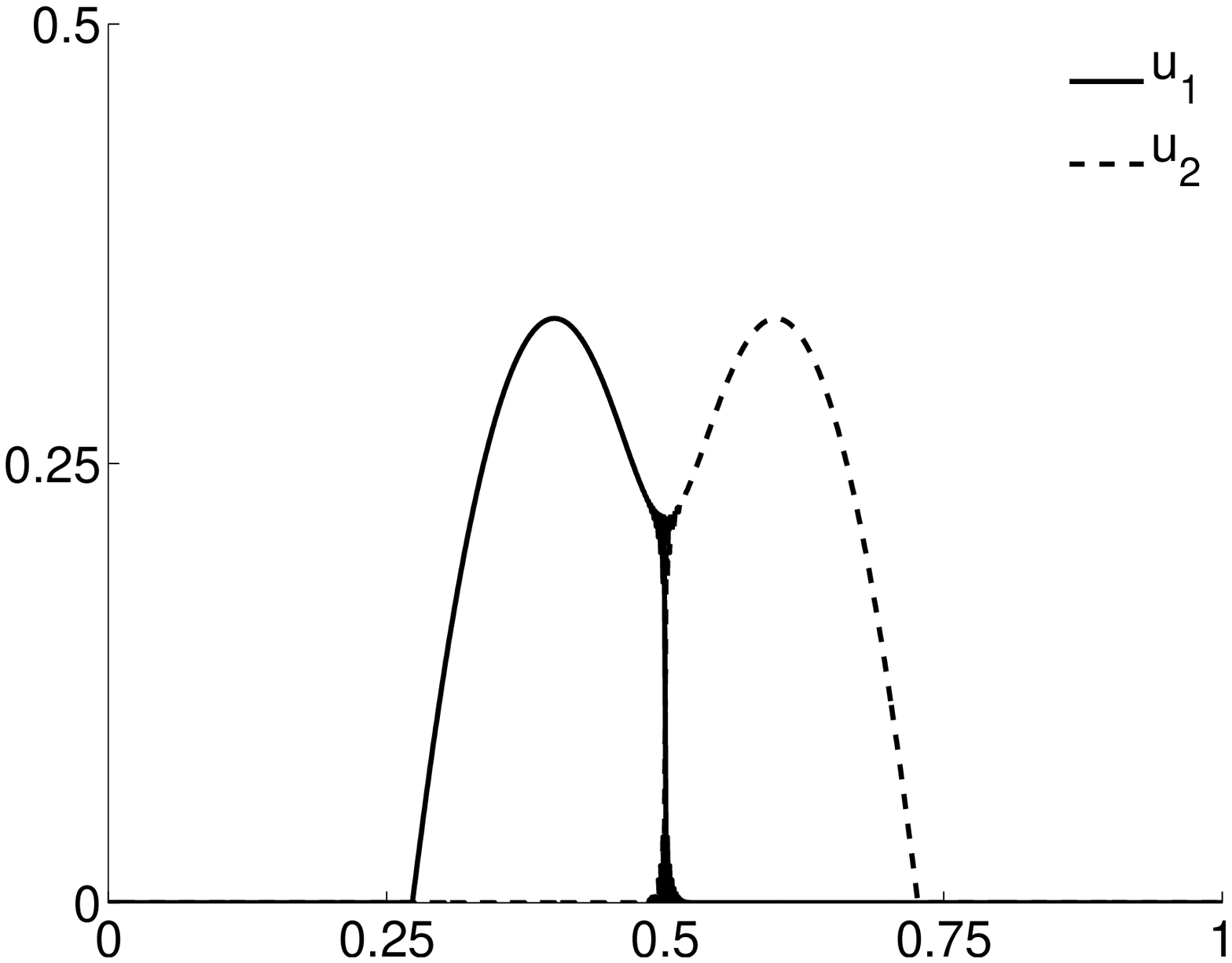}}
 %\hspace{0.5cm}
 \subfigure[$t=0.17$]%
 {\includegraphics[width=4.25cm,height=2.7cm]{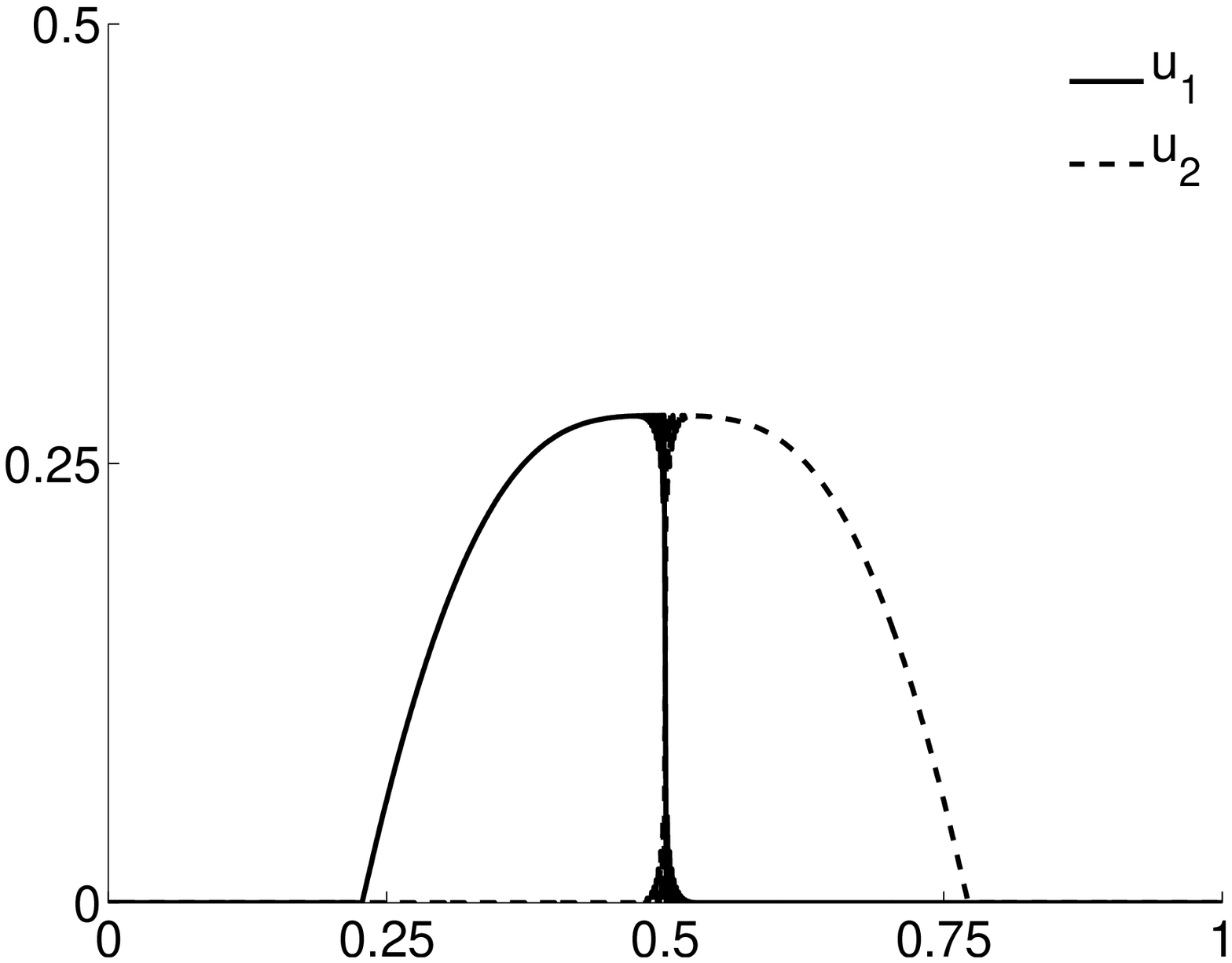}}
  \caption{{\small Experiment 2. The case of  semi-definite positive matrix $(a_{ij})$}} 
\label{exp2_fig}
\end{figure}
\begin{figure}[t]
\centering
 \subfigure[$t=0$]%
 {\includegraphics[width=4.25cm,height=2.7cm]{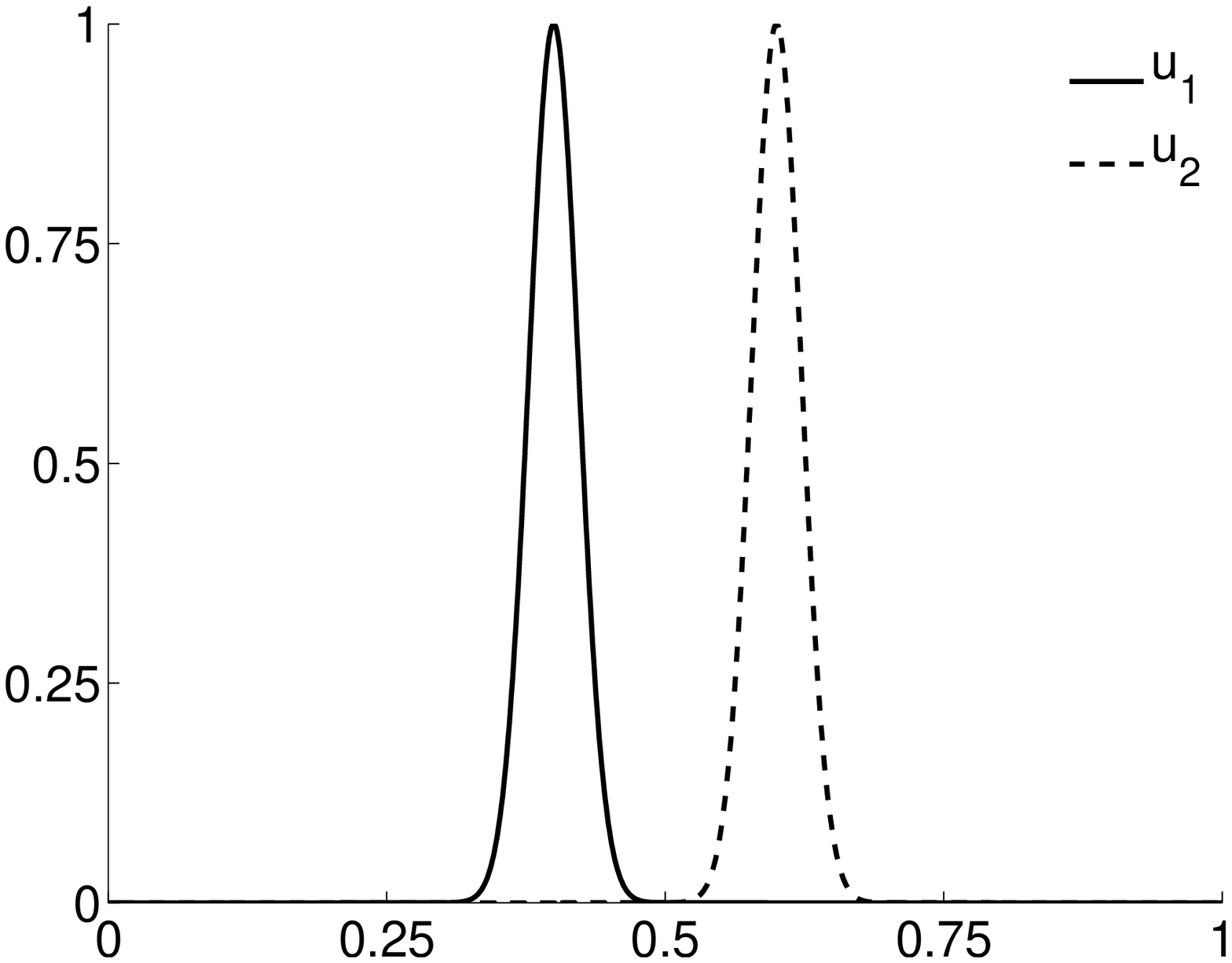}}
%  \hspace{0.5cm}
%  \subfigure[$t=0.01$]
%  {\includegraphics[width=6.25cm,height=4.5cm]{exp3b_fig}} \\
 \subfigure[$t=0.05$]
 {\includegraphics[width=4.25cm,height=2.7cm]{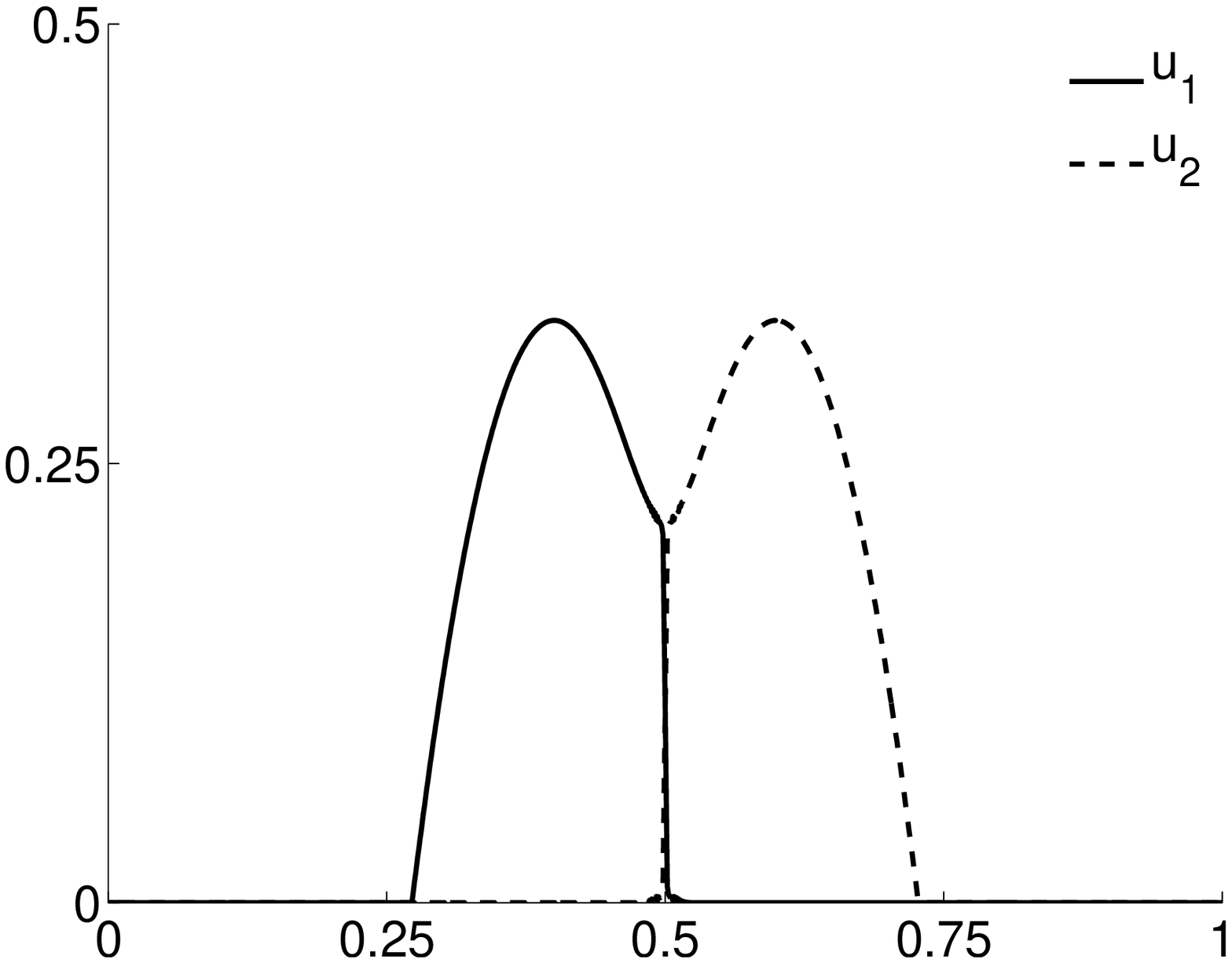}}
%\hspace{0.5cm}
 \subfigure[$t=0.17$]
 {\includegraphics[width=4.25cm,height=2.7cm]{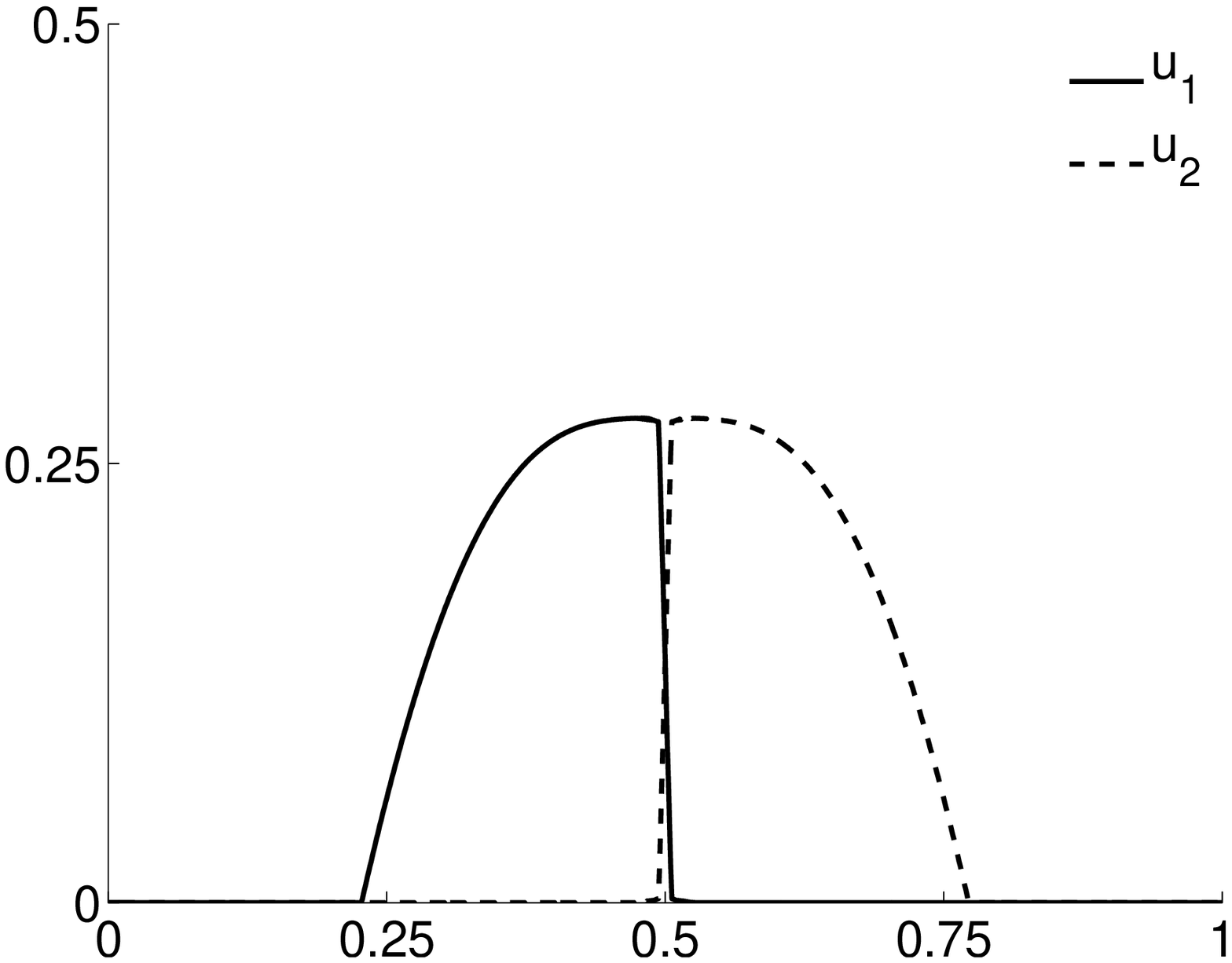}}
  \caption{{\small Experiment 2. Regularizing parameter $\delta=0.001$.}} 
\label{exp3_fig}
\end{figure}
\begin{figure}[ht]
\centering
 \subfigure[$\delta=0$ ]
 {\includegraphics[width=4.25cm,height=2.7cm]{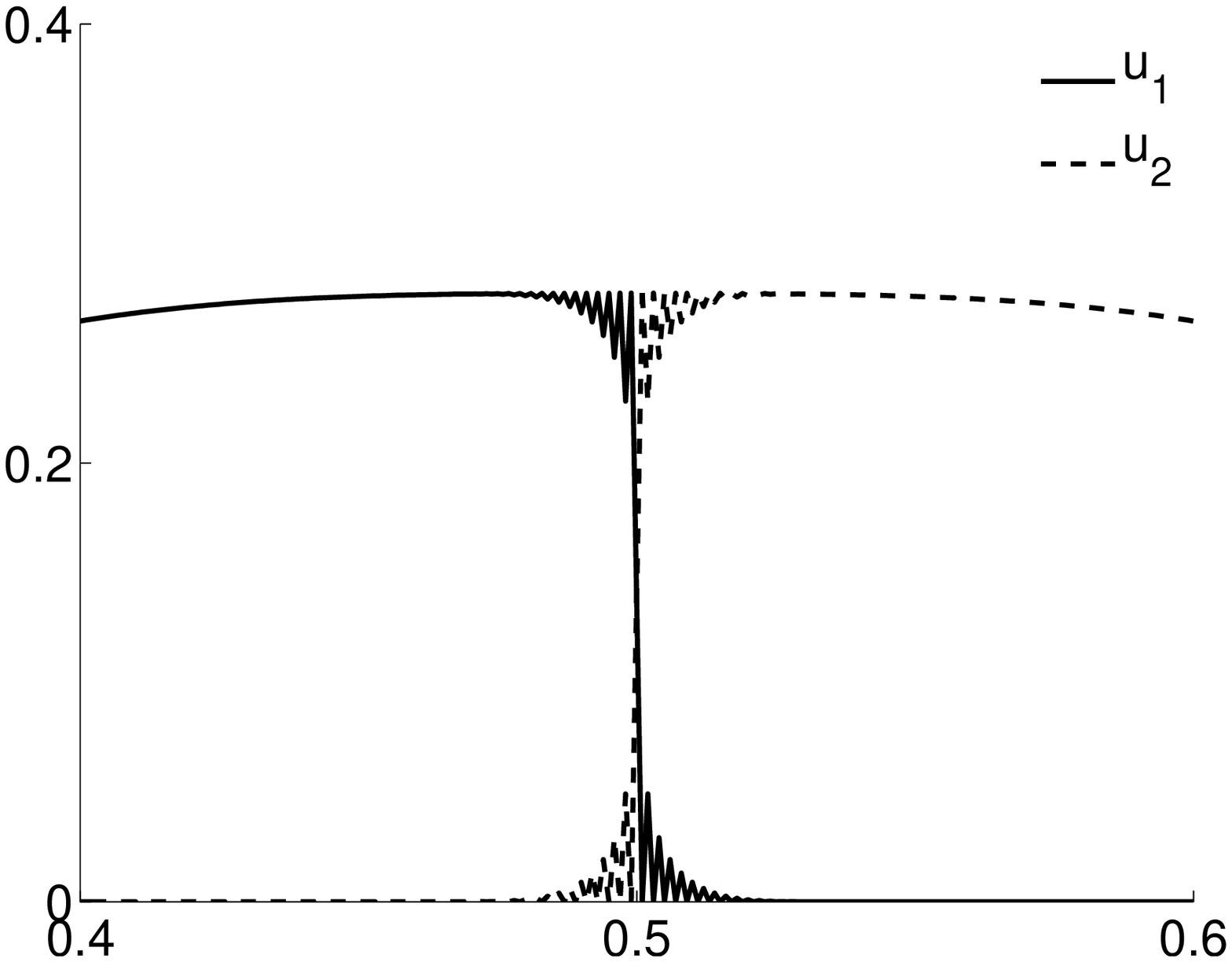}}
 %\hspace{0.5cm}
 \subfigure[$\delta=0.001$]
 {\includegraphics[width=4.25cm,height=2.7cm]{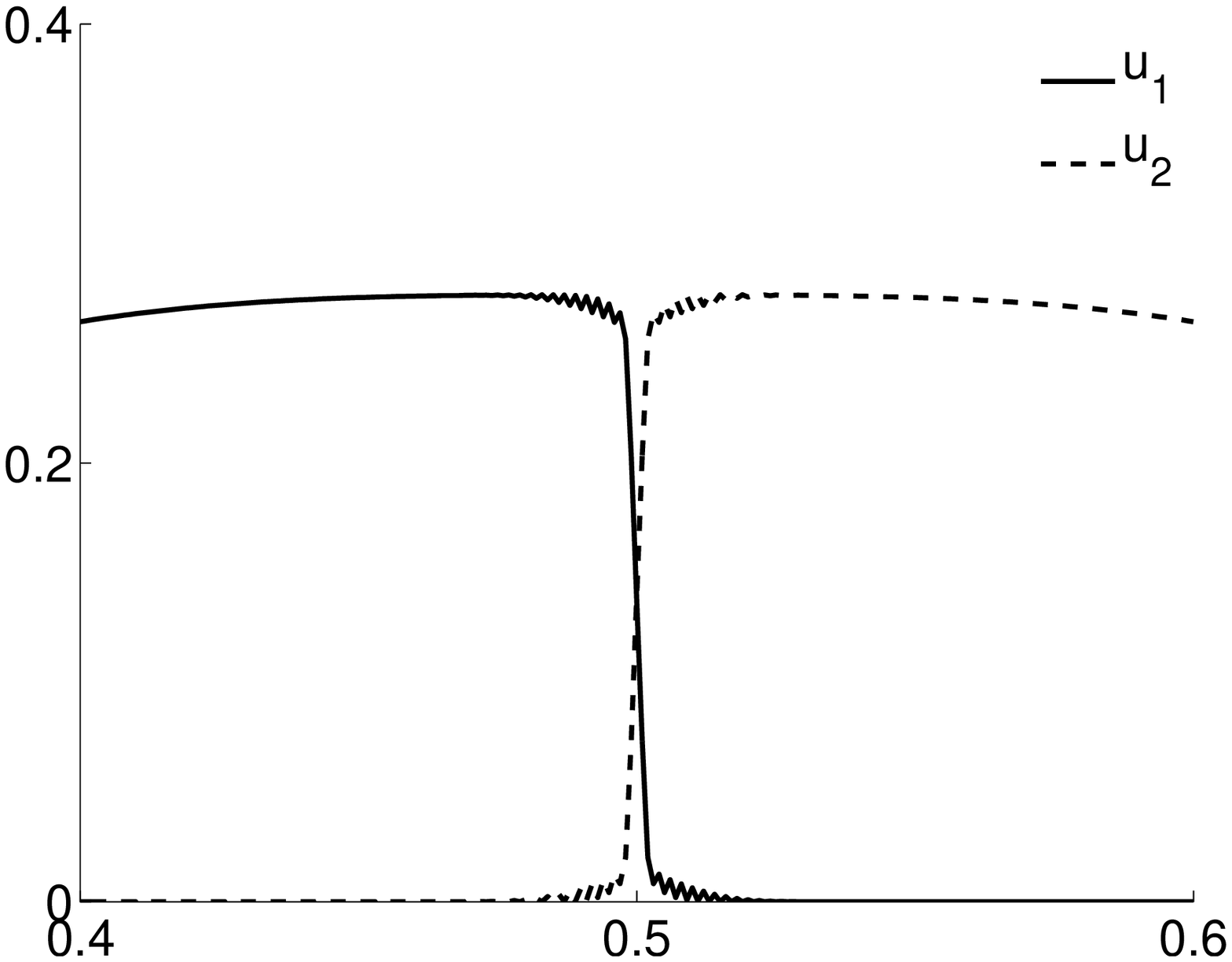}} 
 \subfigure[$\delta=0.01$]
 {\includegraphics[width=4.25cm,height=2.7cm]{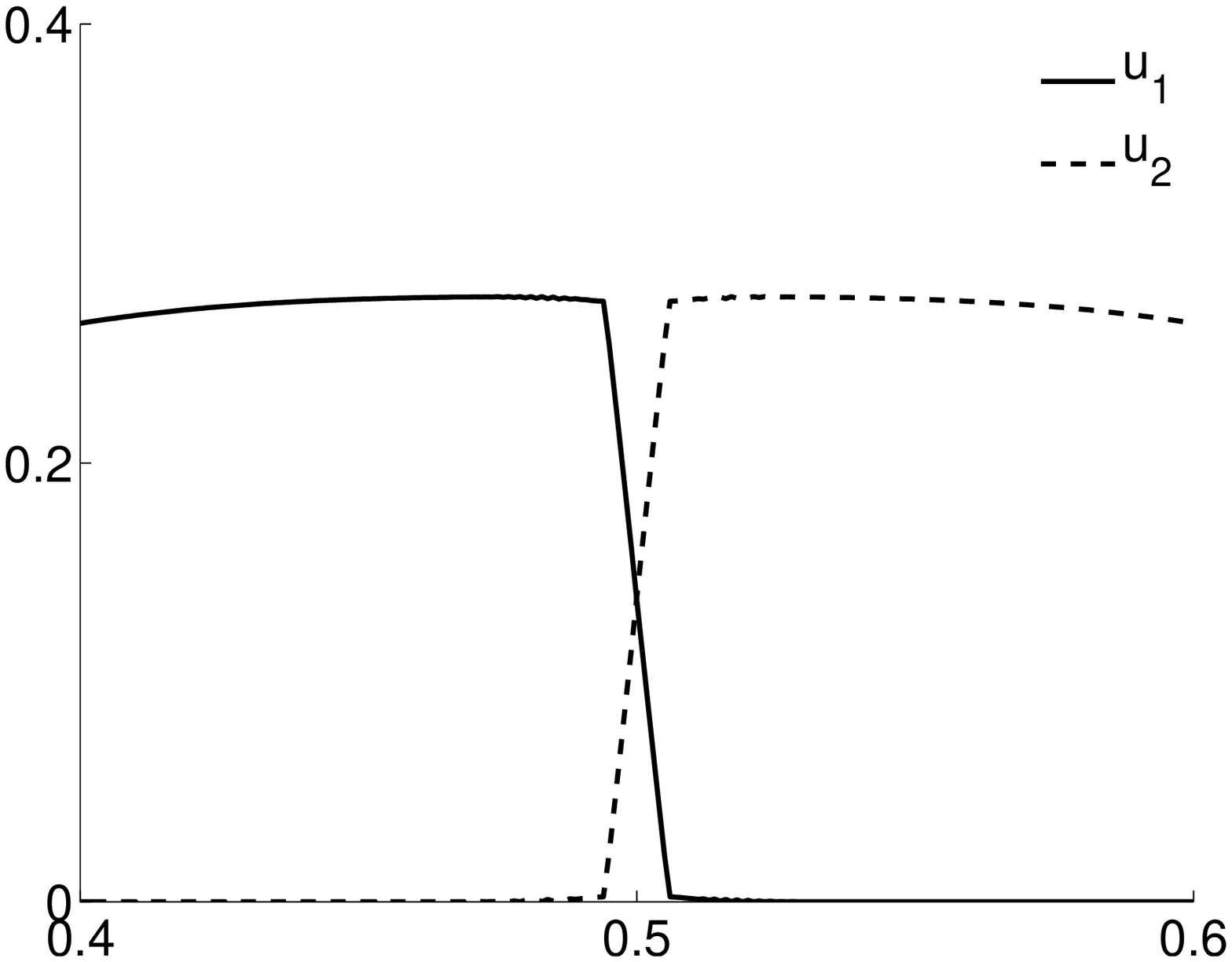}}
%  \hspace{0.5cm}
%  \subfigure[$\delta=0.1$]
%  {\includegraphics[width=6.25cm,height=4.75cm]{exp3e_fig}}
   \caption{{\small Experiment 2. Zoom into the supports intersection region.}} 
\label{exp23_fig}
\end{figure}

%\bigskip 

From the mathematical and numerical point of view, the most interesting situation of problem \fer{eq:pde}-\fer{eq:id} is that of the degenerate case covered by assumptions (H'), i.e. when 
$a_{ij}=a>0$ and, therefore, matrix $(a_{ij})$ is only semi-definite positive. In this particular case, the following
property holds. Let $(u_1^{(\delta)},u_2^{(\delta)})$ be a solution of problem  \fer{eq:pde}-\fer{eq:id} 
with $J_i$ replaced by  $J_i^{BT,\delta}$, see \fer{def.NLflowapp}, to which we refer to as Problem (P)$_\delta$,
and assume (H'). Then $u=u_1^{(\delta)}+u_2^{(\delta)}$ solves
\begin{align}
& \pt u -\Div J^{(\delta)}(u)=f(u) && \qtext{in }Q_T=\O\times(0,T), 
  \label{prob:add.ec}\\
&  J^{(\delta)}(u)\cdot n =0 && \qtext{on }\Gamma_T=\partial\O\times(0,T),
\label{prob:add.bc}\\
& u(\cdot,0)=u_{10}+u_{20} && \qtext{in }\O,	\label{prob:add.ic} 
\end{align}
for $J^{(\delta)}(u)=(a+\delta) u \grad u +b q u +c\grad u$ and 
$f(u)= u (\alpha -\beta u)$, which is a uniformly parabolic problem in view of (H'). 

%\bigskip 

In the following experiments we take, unless otherwise stated, $\O=(0,1)$, $b_i=c_i=0$, and 
$u_{i0}= \exp((x-x_i)^2/0.001)$, $f_i=0$ for $i=1,2$, with $x_1=0.4$ and $x_2=0.6$.
We chose a larger tolerance parameter 
for the fixed point algorithm than in the previous experiment,
$\text{tol}=10^{-4}$, in view of the slow convergence observed for the discretization parameters 
$M=1001$ and $\tau = 10^{-5}$. 
Although the initial data  do not satisfy condition (H')$_4$, this does not seem to affect the convergence or stability of the algorithm for the different cases under study.

%\bigskip 

\no\emph{Experiment 2.} 
We run experiments for solving problem \fer{eq:pde}-\fer{eq:id} with coefficients $a_{ij}=1$, for $i,j=1,2$, and the corresponding regularized version given by Problem (P)$_\delta$.  
We set the final time to $0.17$ and investigate the behavior of solutions during the transient state.

 Due to diffusion, after a while the supports of $u_i$ intersect with each other at one point. In this moment, an important qualitative difference arises between the solutions of the degenerate and the regularized problems. For $\eps=0$, no mixture of populations is observed at subsequent times, and a steep 
gradient or discontinuity is formed at the so-called contact inhibition point. Numerical instabilities are clearly seen
around this point, see Fig.~\ref{exp2_fig}.
However, since $u_1+u_2$ is a solution of problem \fer{prob:add.ec}-\fer{prob:add.ic}   and therefore smooth and non-negative (Barenblatt type, see Theorem~\ref{th:existence_particular}) these instabilities must remain bounded.

In Fig.~\ref{exp3_fig} we plot the solution of problem (P)$_\delta$, for $\delta=0.001$, approximating the solution of the degenerate problem shown in Figure~\ref{exp2_fig}. As it can be seen, instabilities do not arise (at the scale of the plot) for this regularized problem. We also observe that the components of the solution mixes
in an interval of $\delta-$dependent length.

Finally, in Fig.~\ref{exp23_fig} we show a zoom of the solutions of problem \fer{eq:pde}-\fer{eq:id} and problems (P)$_\delta$, for several choices of $\delta$, around the intersection point $x=0.5$. As suggested by estimate  \fer{est:lin} (with $c=0$), the square of the $L^2$ norm of the gradient of the solutions seems to be proportional to, and not just bounded by, $1/\delta$.    

%\bigskip 

\begin{figure}[t]
\centering
 \subfigure[$t=0$ ]
 {\includegraphics[width=4.25cm,height=2.7cm]{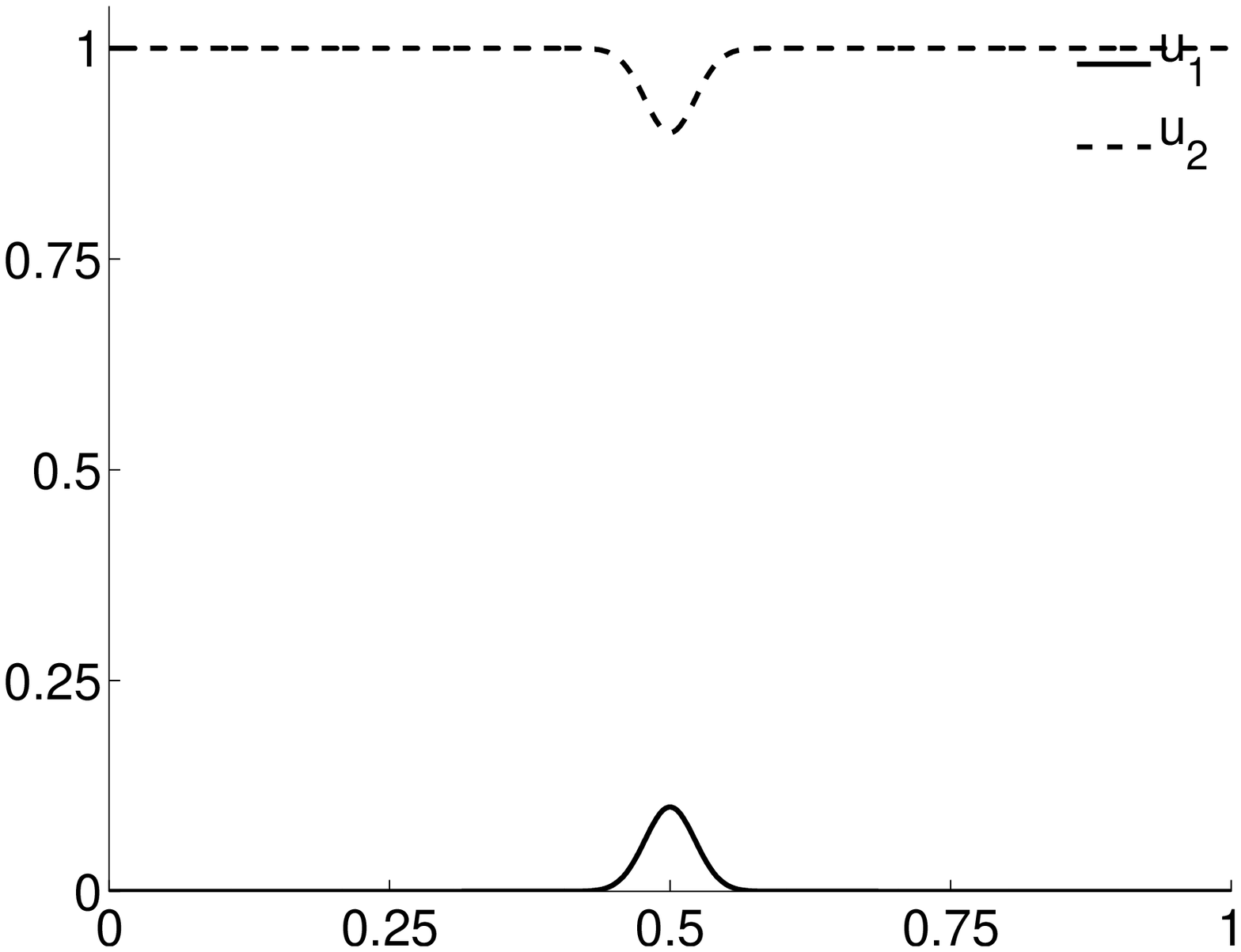}}
 %\hspace{0.5cm}
 \subfigure[$t=5$]
 {\includegraphics[width=4.25cm,height=2.7cm]{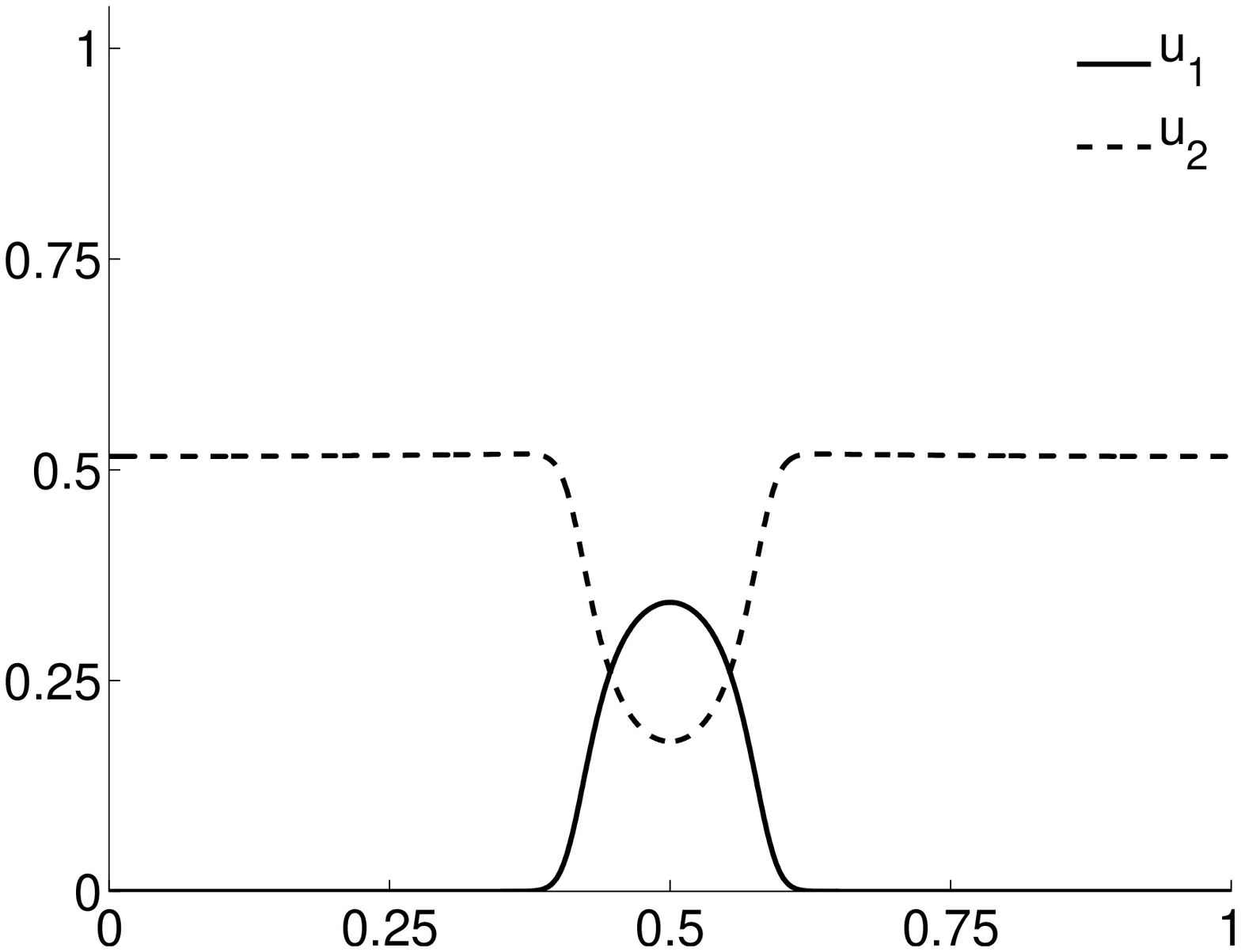}} 
 \subfigure[$t=7$]
 {\includegraphics[width=4.25cm,height=2.7cm]{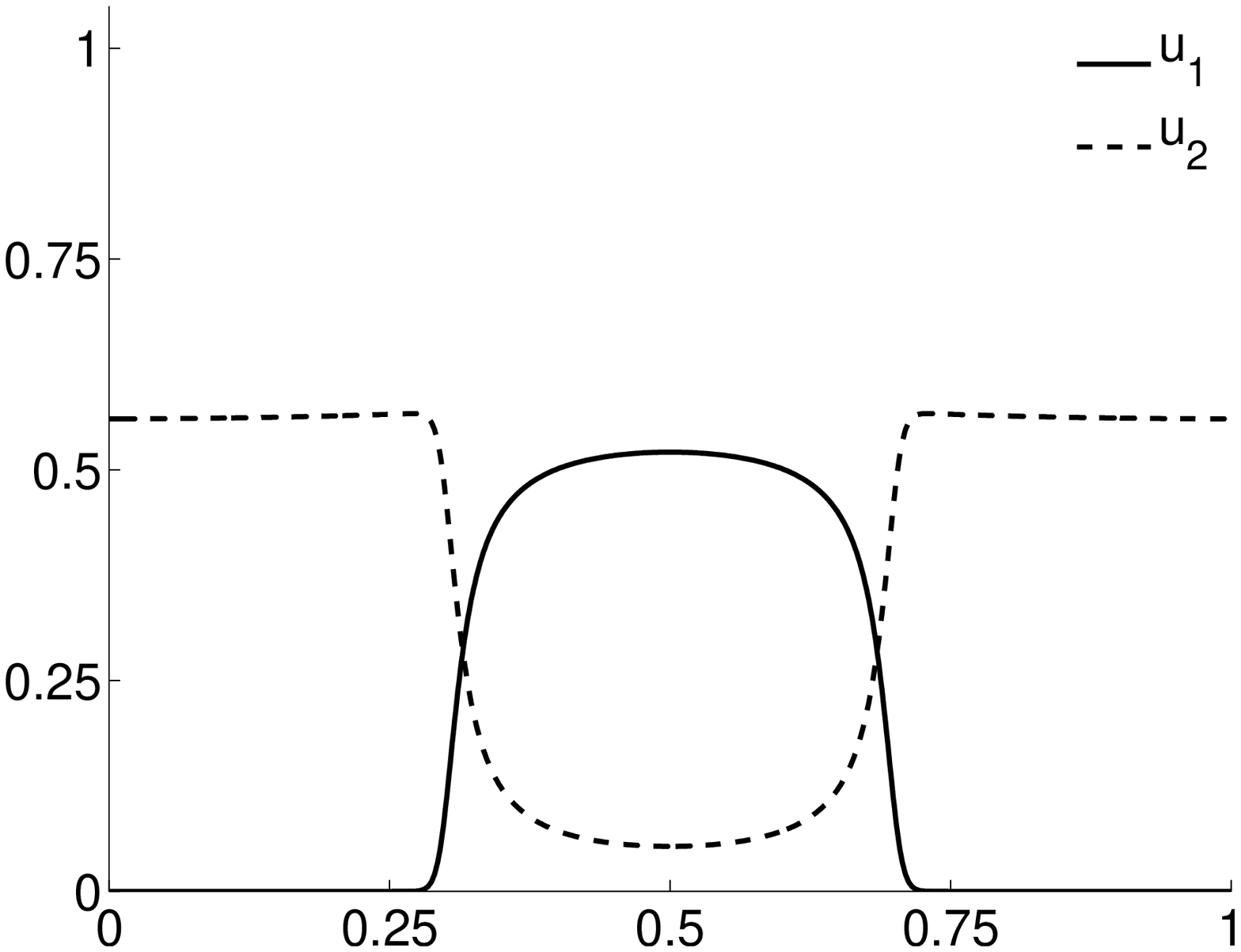}}
%  \hspace{0.5cm}
%  \subfigure[$\delta=0.1$]
%  {\includegraphics[width=6.25cm,height=4.75cm]{exp3e_fig}}
   \caption{{\small Experiment 3. Evolution of the invasion.}} 
\label{exp6_fig}
\end{figure}

\no\emph{Experiment 3.} In this experiment, we look at the question (Q2) stated in \cite{bertsch10}, in which 
the invasion of one population (mutated abnormal cells) over an initially dominant population (normal cell) 
is produced. The initial data is taken as $u_{10}(x)=0.1\exp((x-x_i)^2/0.001)$ and $u_{20}=1-u_{10}$. The 
Lotka-Volterra competitive term is taken of the usual form \fer{def:reaction}, with 
$\alpha_1=\alpha_2=1$, $\beta_{11}=\beta_{12}=1$ and $\beta_{21}=\beta_{22}=2$. In Fig.~\ref{exp789_fig}
we show the initial distributions and two instants of the transient state, in which the pressure exerted by
the mutant population drives the system to a change of equilibrium. The steady state, not shown in the figure, 
is of extinction of the normal cells.

\bigskip 

\no\emph{Experiment 4.} In the following simulations we investigate other parameter ranges out of those stated in (H'). 
We take the same parameters than in Experiment~2 but
\begin{enumerate}
 \item Changing matrix $(a_{ij})$ to $a_{11}=a_{12}=3$, $a_{21}=a_{22}=1$, so still positive semi-definite.
 We may see a transient state of the solution in the left panel of Fig.~\ref{exp789_fig}. 
 
 \item Changing the transport coefficients to $d_1=1$ and $d_2=10$, and setting $q(x)=-3(x-0.5)$. A transient state is plotted in the center panel of Fig.~\ref{exp789_fig}. 
\end{enumerate}
As we see, both set of parameters produce continuous solutions with discontinuous gradients at the contact inhibition point. This, again, suggests that our conditions (H') are not optimal. In particular, a solution to the case 1. may be constructed by using Lagrangian coordinates, see \cite{bertsch85}.

%\newpage

\section{Proofs of the theorems}

\subsection{Proof of Theorem~\ref{th:existence_original}} 
To make the entropy inequality \fer{ineq:entropy} rigorous one has to go through a regularization procedure. We use the approach introduced by Barrett 
and Blowey \cite{barret04}, even though alternative approaches are also possible, see Chen and J\"ungel \cite{chen04}. 
Although the results in \cite{chen04,barret04} can not be directly applied to prove Theorem~\ref{th:existence_original},
we use similar techniques in its proof. For the sake of completeness, we replicate some of the arguments used in \cite{barret04}, showing how 
they adapt to problem \fer{eq:pde}-\fer{eq:id} under assumptions (H).

Let $\eps\in (0,1)$ and consider $F_\eps:\R\to [0,\infty)$ given by
\begin{equation}\label{def:F_eps}
F_{\eps} (s):= 
\begin{cases}
\frac{s^2-\eps ^2}{2\eps} + s (\ln \eps -1) +1 & \qtext{if } s\leq \eps ,
\\[1ex] 
s (\ln s-1)+1 & \qtext{if } \eps\leq s\leq \eps^{-1} , \\[1ex]
\frac{\eps (s^2-\eps^{-2}) }{2} + s (\ln \eps^{-1} -1) +1 &
\qtext{if } \eps^{-1}\leq s .
\end{cases}
\end{equation}
Notice that function $F$ given in \fer{def:entropy} is defined in $[0,\infty)$, whereas $F_\eps $ is defined in the whole real line, $\mathbb{R}$. % when $\eps\in (0,1)$. 
Besides, $F_{\eps}\in\mathcal{C}^{2,1} (\R) $,
so we may as well define the H\"older continuous function
\begin{equation}\label{def:lambda_eps,tildeeps}
\lambda _{\eps}(s):= 1/F''_{\eps} (s) .
\end{equation} 
The corresponding regularized version of problem \fer{eq:pde}-\fer{eq:id} reads: For $i=1,2$, find $u_{\eps i}:Q_T\to\R$ such that
\begin{align}
&\pt u_{\eps i}
-\Div J_{\eps i}(u_{\eps 1},u_{\eps 2})
  =f_{\eps i}(u_{\eps 1},u_{\eps 2})  && \qtext{in }Q_T, 
  \label{eq:pde_eps}\\
& J_{\eps i}(u_{\eps 1},u_{\eps 2}) \cdot n =0 && \qtext{on
}\Gamma_T,
\label{eq:bc_eps}\\
& u_{\eps i}(\cdot,0)=u_{i0} && \qtext{in }\O,	\label{eq:id_eps} 
\end{align}
with regularized flow and competitive Lotka-Volterra functions given by
\begin{align}
& J_{\eps i}(v_1,v_2) = \lambda _{\eps } (v_i)  \big(a_{i1}\grad v_1 +a_{i2}\grad
v_2 +b_i q ) + c_i\grad v_i ,
   \label{def:flow_eps}\\
&  f_{\eps i} (v_1,v_2) =  \alpha_{i}v_i -\big(\beta_{i1} \lambda _{\eps } (v_1)
 + \beta_{i2} \lambda _{\eps } (v_2)\big) \lambda _{\eps }
(v_i) .
\label{def:reaction_eps}
\end{align}

% FULLY DISCRETE PROBLEM
\subsubsection{Finite element approximation }\label{subsec:fem}
We consider a fully discrete approximation using finite elements in space and backward finite differences in time. We consider  a quasi-uniform family of meshes of $\Omega$ (polygonal), $\{\mathcal{T}_h\} _h$, composed by right-angled tetrahedra, with parameter $h$ representing its diameter.
We introduce the finite element space of piecewise $\mathbb{P}_1$-elements:
$$
S^h = \{ \chi\in \mathcal{C}(\overline{\Omega} ) ; \, \chi |_K\in\mathbb{P}_1\,\text{ for all } K\in\mathcal{T}_h \} .
$$
The Lagrange interpolation operator is denoted by $\Pi ^h : \mathcal{C}(\overline{\Omega} ) \to S^h$. We also introduce 
the discrete semi-inner product on $\mathcal{C}(\overline{\Omega} ) $ and its induced discrete seminorm:
$$
(\eta_1,\eta_2)^h= \int_{\Omega} \Pi^h(\eta_1\eta_2) ,\quad |\eta|_h=\sqrt{(\eta,\eta)^h}.
$$
Finally,  $Q^h: L^2(\Omega)\to S^h$ stands for the $L^2(\Omega)$-projection.

For each $\eps\in (0,1)$ we consider the construction of the linear operator $\Lambda _{\eps} : S^h\to L^{\infty}(\Omega)^{m\times m}$  given in \cite{barret04,grun00} which, for all $z^h\in S^h$ and a.e. in $\O$, has a symmetric and positive image $\Lambda _{\eps} z^h$, and satisfies $(\Lambda _{\eps}z^h) \nabla\Pi ^h (F_{\eps}'(z^h)) = \nabla z^h$.
 Then, due to the right angled constraint requirement, the following bound holds 
\begin{equation}\label{props2-L_eps}
|\nabla \Pi ^h \lambda_{\eps} (\chi) |^2_{1} \leq (\nabla\chi , \nabla \Pi ^h \lambda_{\eps} (\chi) ) \quad\text{for all } \chi\in S^h .
\end{equation}

For the time discretization, we take a possibly non-uniform partition of $[0,T]$ in $N$ subintervals: $0=t_0<t_1<\ldots <t_{N-1}<t_N=T$. We denote the time steps by $\tau _n=t_n-t_{n-1}$ ($n=1,\ldots , N$), and $\tau = \max _{n=1,\cdots,N} \tau _n $.

For the discrete problem we need more regularity on the coefficients than that assumed in (H). Therefore,  
we introduce sequences of nonnegative functions 
$a_{ij}^\sigma,~c_i^\sigma,~\alpha_i^\sigma,~\beta_{ij}^\sigma \in \mathbb{P}_1([0,T];\mathcal{C}(\bar \O))$,
as well as functions $b_i^ \sigma \in \mathbb{P}_1([0,T];\mathcal{C}(\bar \O))$ 
and  $q ^\sigma \in ( \mathbb{P}_1([0,T];\mathcal{C}(\bar \O)))^m $
for $\sigma >0 $, such that, as $\sigma\to 0$,
\[
\begin{array}{l}
a_{ij}^\sigma \to a_{ij},\quad b_i^\sigma \to b_i ,\quad c_i^\sigma \to c_i ,\quad \alpha_i^\sigma \to \alpha_i ,
\quad \beta_{ij}^\sigma \to \beta_{ij}\qtext{strongly in }L^\infty(Q_T), \\[.2em]
q ^\sigma \to q \qtext{strongly in } (L^2(Q_T))^m,
\end{array}
\]
and satisfying \fer{H:def_pos} uniformly in $\sigma$ (i.e. with $a_0$ a constant independent of $\sigma$). We use the following notation for the time-space discretization of coefficients:
\[
a^{\sigma n}_{ij}=\Pi^h (a_{ij}^\sigma (\cdot,t_n)),\quad q^{\sigma n} =(\Pi^h (q^\sigma _1(\cdot,t_n)),\ldots,\Pi^h(q^\sigma _m(\cdot,t_n))),\qtext{etc.}
\]

Finally, we collect here some restrictions on the discretization-regularization parameters that we shall use in this section:
\begin{equation}
 \label{hd:param}
 \eps\in(0,{\rm e}^{-2} ), \quad h>0,\quad \sigma>0,\quad \omega \tau \leq 1-\delta < 1,\qtext{for some }\delta>0
\end{equation}
with $
\omega = \max _{i=1,2} \{ 2\nor{\alpha _i}_{L^\infty(Q_T)} + 
\nor{\beta_{i1}}_{L^\infty(Q_T)} + \nor{\beta_{i2}}_{L^\infty(Q_T)}\}$.

\subsubsection{The discrete problem}

In this subsection we prove the existence of solutions of the fully discrete problem 
corresponding to problem \fer{eq:pde_eps}-\fer{eq:id_eps} and deduce uniform estimates on the solutions 
which will allow us to pass to the limit 
in the discretization-regularization problem to obtain a solution of the continuous problem.
Along this subsection we omit the superindex $\sigma$ in  the unknowns for clarity in the notation. 
\begin{lemma}\label{th:existence-step}
Assume \fer{hd:param} and let, for $n\geq1$, $u^{n-1}_{\eps i}\in S^h$, being  $u^0_{\eps i}:=Q ^h u_{i0}$.  Then, there exists $u^{n}_{\eps i}\in S^h$ solution of the $n$-th step of the fully discrete problem 
\begin{equation}\label{eq:pde_discr}
\begin{array}{l}
\big( \dfrac{u^n_{\eps i}-u^{n-1}_{\eps i}}{\tau _n} , \chi )^h
+ ( \Lambda _{\eps } (u^n_{\eps i})  \big(a_{i1}^{\sigma n} \grad u^n_{\eps 1}
+a_{i2}^{\sigma n}\grad u^n_{\eps 2} +b_i^{\sigma n} q^{\sigma n} \big) + c_i^{\sigma n}\grad u^n_{\eps i}
, \grad \chi )^h
  = \\[2ex]
\hspace*{1cm} = (\alpha_{i}^{\sigma n} u^n_{\eps i} - \lambda _{\eps } (u^n_{\eps i})
\big( \beta_{i1}^{\sigma n} \lambda _{\eps } (u^{n-1}_{\eps 1}) + \beta_{i2}^{\sigma n} \lambda
_{\eps } (u^{n-1}_{\eps 2}) \big) , \chi )^h , %\qquad\forall \chi\in S^h .
\end{array}
\end{equation}
for every $ \chi\in S^h $, and satisfying, for a constant $C$ independent of $\eps$, $h$,  $\tau$ and $\sigma$,
$$
\begin{array}{l}
\displaystyle\max_{n=1,\ldots,N} \Big( \sum_{i=1}^{2} (F_{\eps}
(u^n_{\eps i}),1)^h +
\eps^ {-1} |\Pi^h [u^n_{\eps i}]_{-}|^2_0 +  |u^n_{\eps i}|_{0,1} \Big)
+ a_0 \sum _{n=1}^{N} \tau _n \sum _{i=1}^{2} ||u^n_{\eps i}||_1^2
\leq C,
\end{array}
$$
and, for  $r=2(m+1)/m$ and $p=2(m+1)/(2m+1)$, 
\begin{align*}
\sum _{n=1}^{N} \tau _n \sum _{i=1}^{2} \displaystyle\Big( |\Lambda _\eps
(u^n_{\eps i})|^r_{0,r}
 + |\Pi^h(\lambda _{\epsilon } (u^n_{\eps i}))|^r_{0,r}
 + |u^n_{\eps,i}|^r_{0,r}  
 +  |\lambda _{\eps } (u^n_{\eps i})|^r_{0,r}  \\
+ ||\mathcal{G}
(\frac{u^n_{\eps i}-u^{n-1}_{\eps i}}{\tau_n})||^p_{1,p} \displaystyle\Big)
\leq C .
\end{align*}
\end{lemma}

\no\emph{Proof of Lemma \ref{th:existence-step}.}
We split the proof into three steps. 

\noindent\textbf{Step 1. }We prove the existence of solutions of the discrete problem with
a proof by contradiction. Let us define $\mathcal{A}\equiv (A_1,A_2) : S^h\times S^h\to S^h\times S^h$  by
\begin{equation}\label{proof:existence-step_0}
\begin{array}{l}
(A_i (v_1,v_2), \chi )^h = 
( v_{i} - u^{n-1}_{\eps i} , \chi )^h
+ \tau _n \, ( \Lambda _{\eps } (v_{i})  \big(a_{i1}^{\sigma n} \grad v_{1} + a_{i2}^{\sigma n}\grad v_{2} + b_i^{\sigma n} q^{\sigma n} \big) 
, \grad\chi )^ h \\[2ex]
\hspace*{1cm} + \tau_n ( c_i^{\sigma n}\grad v_{i} , \grad \chi )^h 
 - \tau _n\, (\alpha_{i}^{\sigma n} v_{i} - \lambda _{\eps } (v_{i}) \big( \beta_{i1}^{\sigma n} \lambda _{\eps } (u^{n-1}_{\eps 1}) + \beta_{i2}^{\sigma n} \lambda _{\eps } (u^{n-1}_{\eps 2}) \big) , \chi )^h , %\qquad\forall \chi\in S^h .
\end{array}
\end{equation}
for every $\chi\in S^h$. Then, the $n$-th step of the fully discrete problem, (\ref{eq:pde_discr}), consists of finding $u^n_{\eps i}\in S^h$ ($i=1,2$) such that
\begin{equation}\label{proof:existence-step_1}
\mathcal{A} (u^n_{\eps 1},u^n_{\eps 2}) = (0,0) .
\end{equation}
Suppose a solution does not exist and let $R>0$ be such that 
$$
\sum_{i=1}^2 |A_i (v_1,v_2)| > 0 \qtext{for all } (v_1,v_2)\in S^h_R = \{ (v_1,v_2)\in S^h\times S^h ; \, |v_1|_h^2+|v_2|_h^2\leq R^2 \} .
$$
Consider the function 
$B=(B_1,B_2): S^h_R\to S^h_R$ given by 
$$
B_i (v_1,v_2) := -R \, A_i (v_1,v_2) \big(\sum _{j=1}^2 |A_j (v_1,v_2)|_h^2\big)^ {-1/2}.
$$
We have: (i) $S_R^ h$ is a convex a compact subset of $S^h\times S^h$, (ii) $B$ is continuous in $S_R^ h$,
since $\mathcal{A}|_{S_R^ h}$ is well defined and continuous, and (iii) $B(S_R^h)\subset S_R^ h$. Then, Brouwer's fixed-point theorem guarantees the existence of $(w_1,w_2)\in S^h_R$ such that $B(w_1,w_2)=(w_1,w_2)$ which, in particular, satisfies $|w_1|_h^2+|w_2|_h^2=R^2$. Taking $v_1=w_1$, $v_2=w_2$ and $\chi=\Pi^h F_{\eps}'(w_i)$ in (\ref{proof:existence-step_0}), we obtain, using our assumption \fer{hd:param} on $\tau$,
\begin{equation}
 \label{contradiction}
\sum _{i=1}^2 (A_i (w_1,w_2), \Pi^h F_{\eps}'(w_i) )^h \geq \frac{\eps}{4} R^2 - C , 
\end{equation}
with $C$ a constant independent of $\eps$, $R$, $w_1$ and $w_2$. Then,  for $R>0$ large enough,  the following contradiction arises:
On one hand, by \fer{contradiction} and using that $(w_1,w_2)$ is a fixed point of $\mathcal{B}$, we obtain
$$
\sum _{i=1}^2 (w_i,F_{\eps}'(w_i))^h \leq -R \, \big(\frac{\eps}{4} R^2 - C\big)\big(\sum _{j=1}^2 |A_j (w_1,w_2)|_h^2\big)^ {-1/2} <0.
$$%th:existence-step
On the other hand, by standard properties of function $F_\eps$, we deduce 
$$
\sum _{i=1}^2 (w_i,F_{\eps}'(w_i))^h \geq \sum _{i=1}^2 \left( (F_{\eps}(w_i)-F_{\eps}(0),1)^h +\frac{\eps}{2} \, |w_i|_h^2 \right) \geq -2\abs{\O} + \frac{\eps}{2} \, R^2 >0.
$$

\noindent\textbf{Step 2. }We now pass to the proof of the first estimate of Lemma~\ref{th:existence-step}.
Taking  $\chi=\Pi^hF_\eps '(u^ n_{\eps i})$ in (\ref{eq:pde_discr}) and summing over $i=1,2$ we deduce 
\[
\displaystyle\sum_{i=1}^{2} (1-\omega  \tau _n )\, (F_{\eps} (u^n_{\eps i}),1)^h 
+ \frac{a_0\tau _n}{ 2} \sum _{i=1}^{2} |u^n_{\eps i}|_1^2 
\leq
(1+\tau_n ) \,\displaystyle\sum_{i=1}^{2} (F_{\eps} (u^{n-1}_{\eps i}),1)^h  + C\, \tau _n ,
\]
from where we obtain 
$$
(1-\omega \tau_n) \sum_{i=1}^2 (F_{\eps} (u^n_{\eps i}), 1)^h 
\leq 
\sum_{j=1}^{n-1} \sum_{i=1}^2 (\tau_{j+1}+\omega  \tau_{j}) (F_{\eps} (u^j_{\eps i}), 1)^h 
+ \sum_{i=1}^2 (F_{\eps} (u^0_{\eps i}), 1)^h + T \, C .
$$
By \fer{hd:param} and the discrete Gronwall's lemma, we get
\begin{equation}
 \label{mm1}
\max_{n=1,\ldots,N} \sum_{i=1}^2 (F_{\eps} (u^n_{\eps i}), 1)^h 
\leq 
\frac{1}{\delta} \,
\big(1+ T
%\frac{T}{\delta}
(1+\omega ) \mathrm{e}^{(1+\omega )\,T %/\delta
}\big) 
\big(\sum_{i=1}^2 (F_{\eps} (u^0_{\eps i}), 1)^h \, (\tau +1)
+ T \, C \big) .%\leq C .
\end{equation}
Similarly, choosing $\chi=1$ as test function in (\ref{eq:pde_discr}), 
leads to
\begin{equation}
 \label{mm2}
 \max_{n=1,\ldots,N}(u^n_{\eps i} ,1)^h \leq \frac{1}{\delta} \big( 1+\frac{\omega T}{\delta}\mathrm{e}^{\frac{\omega T}{\delta}} \big) \, (u^0_{\eps i} ,1)^h .
\end{equation}
Since 
$|u^n_{\eps i}|_{0,1} \leq (u^n_{\eps i} ,1)^h + 2 \, |\Pi^h [u^n_{\eps i}]_{-}|_{0,1}$, and 
$|\Pi^h [u^n_{\eps i}]_{-}|_{0,1} \leq C\,\eps^{1/2} \, (F_{\eps} (u^n_{\eps i}),1)^h$, 
we also obtain, using \fer{mm1}, \fer{mm2} and standard properties of function $F_\eps$,
\begin{equation}
 \label{mm3}
\max_{n=1,\ldots,N} |u^n_{\eps i}|_{0,1} + \eps^ {-1} |\Pi^h [u^n_{\eps i}]_{-}|_{0,1}^2 \leq C.
\end{equation}
 
\noindent\textbf{Step 3. }We finish the proof of the lemma proving the last estimate.
Using the properties of the mapping of the reference element onto an element $K\in\mathcal{T}_h$ as well as   Sobolev embedding theorem (for $r=2 (m+1)/m$), we obtain
$$
\begin{array}{l}
|\Lambda_{\eps}(u^n_{\eps i})|_{0,r}^r\leq C\displaystyle\sum_{K\in\mathcal{T}_h} \int_K 1\,\, |\Pi^h\lambda_{\eps}(u^n_{\eps i})|^r_{0,\infty ,K} \leq C \, |\Pi^h\lambda_{\eps}(u^n_{\eps i})|^r_{0,r ,\Omega } \\[1ex]
\hspace*{3cm} \leq C\, |\Pi^h\lambda_{\eps}(u^n_{\eps i})|^{r-2}_{0,1} \, ||\Pi^h\lambda_{\eps}(u^n_{\eps i})|^2_1 .
\end{array}
$$
By Poincar\'e inequality and (\ref{props2-L_eps}),
$
||\Pi^h\lambda_{\eps}(u^n_{\eps i})||^{2}_{1}\leq C \, (|\Pi^h\lambda_{\eps}(u^n_{\eps i})|^{2}_{0,1}+ |u^n_{\eps i}|^{2}_{1}).
$
Besides,
$$
|\Pi^h\lambda_{\eps}(u^n_{\eps i})|_{0,1} \leq \eps \int_{\Omega}1 + |u^n_{\eps i}|_{0,1} .
$$
Therefore,
$
|\Lambda_{\eps}(u^n_{\eps i})|_{0,r}^r\leq C \, (1+|u^n_{\eps i}|^{2}_{1}).
$
Moreover,
$$
|\Pi ^h \lambda_{\eps}(u^n_{\eps i})|_{0,r}^r + | \lambda_{\eps}(u^n_{\eps i})|_{0,r}^r \leq C \, (1+ | u^n_{\eps i}|^{r}_{0,r} ) .
$$
Using  $| u^n_{\eps i}|^{r}_{0,r} \leq C\, || u^n_{\eps i} ||^{2}_{1} $, we deduce 
$$
|\Pi ^h \lambda_{\eps}(u^n_{\eps i})|_{0,r}^r + | \lambda_{\eps}(u^n_{\eps i})|_{0,r}^r \leq C \, (1+ || u^n_{\eps i}|| ^{2}_{1} ) .
$$
Finally, let  $\mathcal{G}: (W^{1,p'} (\Omega))'\to W^{1,p} (\Omega)$ be given by
$$
\int _{\Omega} (\nabla\mathcal{G} v \cdot \nabla w + \mathcal{G} v \, w ) = \langle v , w \rangle %_{W^{1,p} (\Omega)'\times W^{1,p} (\Omega)} 
\quad\text{for all } w\in W^{1,p'} (\Omega)  ,
$$
being $\langle\cdot ,\cdot \rangle $ the duality product ${(W^{1,p'} (\Omega))'\times W^{1,p'} (\Omega)}$.
%Notice that $||\mathcal{G}\cdot ||_{1,p'}$ defines a norm on $(W^{1,p} (\Omega))'$.
From  problem (\ref{eq:pde_discr}) we deduce
$$
\begin{array}{l}
\displaystyle\int _{\Omega} (\nabla\mathcal{G} \frac{u^n_{\eps i}-u^{n-1}_{\eps i}}{\tau_n} \cdot \nabla w + \mathcal{G} \frac{u^n_{\eps i}-u^{n-1}_{\eps i}}{\tau_n} \, w ) = \\[1.5ex]
\hspace*{1cm}
= - ( \Lambda _{\eps } (u^n_{\eps i})  \big(a_{i1}^{\sigma n} \grad u^n_{\eps 1}
+a_{i2}^{\sigma n}\grad u^n_{\eps 2} +b_i^{\sigma n}q^{\sigma n} \big) + c_i^{\sigma n}\grad u^n_{\eps i}
, \grad Q^hw )^h\\[1.5ex]
\hspace*{2cm} +  (\alpha_{i}^{\sigma n} u^n_{\eps i} - \lambda _{\eps } (u^n_{\eps i})
\big( \beta_{i1}^{\sigma n} \lambda _{\eps } (u^{n-1}_{\eps 1}) + \beta_{i2}^{\sigma n} \lambda
_{\eps } (u^n_{\eps 2}) \big) , Q^hw )^h 
\end{array}
$$
for $w\in W^{1,p'}(\Omega )$. In consequence,
$$
\begin{array}{l}
\displaystyle\int _{\Omega} (\nabla\mathcal{G} \frac{u^n_{\eps i}-u^{n-1}_{\eps i}}{\tau_n} \cdot \nabla w + \mathcal{G} \frac{u^n_{\eps i}-u^{n-1}_{\eps i}}{\tau_n} \, w ) \leq \\[1.5ex]
\hspace*{0.5cm}
 \leq C \, ||w||_{1,p'} \Big( ( 1+ |\Lambda _{\eps } (u^n_{\eps i})|_{0,r}) \, (1+\displaystyle\sum_{j=1}^2 |u^n_{\eps j}|_{1}) + \sum_{j=1}^{2} |\lambda _{\eps } (u^n_{\eps i})|_{0,2} \, | \lambda _{\eps } (u^n_{\eps j}) |_{0,r} \Big)
\end{array}
$$
for $w\in W^{1,p'}(\Omega )$, and therefore
$$
%\begin{array}{l}
\displaystyle ||\mathcal{G} \frac{u^n_{\eps i}-u^{n-1}_{\eps i}}{\tau_n}||_{1,p}^p  %\leq \\[1.5ex]
%\hspace*{1cm}
 \leq C \, \Big( \!( 1\!+\! |\Lambda _{\eps } (u^n_{\eps i})|_{0,r}) \, (1\!+\!\displaystyle\sum_{j=1}^2 |u^n_{\eps j}|_{1}) \!+\! \sum_{j=1}^{2} |\lambda _{\eps } (u^n_{\eps i})|_{0,2} | \lambda _{\eps } (u^n_{\eps j}) |_{0,r}\! \Big). 
%\end{array}
$$
The statement follows recalling that $| \lambda_{\eps}(u^n_{\eps i})|_{0,2} \leq C \, (1+ || u^n_{\eps i}|| ^{2}_{1} ) $.
$\Box$

\subsubsection{Passing to the limits} In this subsection we construct the solution to the continuous problem.
% notation 
We make now explicit the dependence of the solution on parameter $\sigma$.
For each $n=1,2,\ldots,N$, we define
\begin{equation}\label{def:ui-t}
u_{\eps   i}^\sigma(t) := \frac{t-t_{n-1}}{\tau _n} u_{\eps   i}^{\sigma n} + \frac{t_n-t}{\tau
_n} u_{\eps   i}^{\sigma (n-1)} \quad\forall t\in [t_{n-1},t_n ] ,
\end{equation}
and also consider
\begin{equation}\label{def:ui+-}
(u_{\eps  i}^\sigma)^{-} (t) := u_{\eps  i}^{\sigma (n-1)} \quad\text{and}\quad (u_{\eps  i}^\sigma)^{+}
(t) := u_{\eps  i}^{\sigma n}\quad\forall t\in (t_{n-1},t_n ] .
\end{equation}
In terms of this notation, the fully discrete problem (which has a solution ensured by Lemma~\ref{th:existence-step}), is written as
\begin{equation}\label{eq:pde_discr-all}
\begin{array}{l}
\displaystyle\int _0^T\Big( ( \pt u_{\eps i}^\sigma, \chi )^h
+ \big( \Lambda _{\eps } ((u_{\eps i}^\sigma)^+)  \big(a_{i1}^\sigma\grad (u_{\eps 1}^\sigma)^+ +a_{i2}^\sigma
\grad (u_{\eps 2}^\sigma)^+ +b_i^\sigma q ^\sigma \big) + c_i^\sigma\grad (u_{\eps i}^\sigma)^+ , \grad \chi \big)^h
\Big)
  = \\[1ex]
\hspace*{1cm} = 
\displaystyle\int _0^T (\alpha_{i}^\sigma (u_{\eps i}^\sigma)^+ - \lambda _{\eps }
((u_{\eps i}^\sigma)^+) \big( \beta_{i1}^\sigma \lambda _{\eps } ((u_{\eps 1}^\sigma)^-) + \beta_{i2}^\sigma
\lambda _{\eps } ((u_{\eps 2}^ \sigma)^-) \big) , \chi )^h \, ,
\end{array}
\end{equation}
$u_{\eps i} \in \mathcal{C} ([0,T]; S^h )$, $i=1,2$, and  for every $ \chi\in L^2((0,T);S^h)$, and satisfying a discrete version of the initial condition:
\begin{equation} \label{eq:id_discr-all}
u_{\eps i}^ \sigma (0) = u_{\eps i}^{0} \in S^h .
\end{equation}
Theorem~\ref{th:existence_original} is a direct consequence of the uniform convergence properties of the 
sequence constructed through the solutions to the fully discrete problem, and stated in the following lemma.
The proof of Lemma~\ref{lem:conv} mimics that of Lemma~3.1 and Theorem~3.1 of \cite{barret04}, and therefore we omit it.
\begin{lemma}\label{lem:conv}
Assume \fer{hd:param} and let %\textbf{donde aparece?} la conv fuerte permite deducir que la solución es >=0
$s\in [2,\infty ]$ if $m=1$, $s\in [2,\infty )$ if $m=2$, and $s\in [2,6)$ if $m=3$.  
Consider regularization and discretization parameters satisfying
\begin{equation}\label{hip:conv}
\sigma\to 0 , \quad \tau\to 0, \quad\text{and}\quad \eps h^{-\frac{m}{m+1}}\to 0 \quad\text{as }
h\to 0 ,
\end{equation}
and the first time step satisfying $\tau _1 \leq C h^2$.
Then, there exist non-negative functions $u_i$, $i=1,2$, with
\begin{equation*}
u_i\in L^2 ( (0,T); H^1 (\Omega )) \cap L^r(Q_T) \cap W^{1,p} ((0,T);( W^{1,p'} (\O))') ,
\end{equation*}
such that any sequence of solutions $u_{\eps i}^\sigma \in \mathcal{C} ([0,T]; S^h )$ ($i=1,2$) of (\ref{eq:pde_discr-all})-(\ref{eq:id_discr-all}) has a subsequence (not relabeled) such that
%(that we identify with the whole sequence) with
\begin{equation*}
\begin{array}{ll}
(u_{\eps i}^\sigma)^{(\pm)} \wto  u_i & \qtext{weakly in } L^2(0,T;H^1(\O))\cap L^r(Q_T),\\[1ex]
\mathcal{G}(\pt u_{\eps i}^\sigma) \wto  \mathcal{G} (\pt u_i) & \qtext{weakly in } L^p(0,T; W^{1,p}(\O)),\\[1ex]
(u_{\eps i}^\sigma)^{(\pm)} \to u_i & \qtext{strongly %and a.e.
 in } L^2((0,T);L^{s}(\O)),\\[1ex]
\lambda_{\eps} ((u_{\eps i}^\sigma)^{\pm}) \to u_i & \qtext{strongly % and a.e. 
in } L^2((0,T);L^{r}(\O)),\\[1ex]
\Pi^h (\lambda_{\eps} ((u_{\eps i}^\sigma)^{\pm})) \to u_i & \qtext{strongly % and a.e.
 in }L^2((0,T); L^{r}(\O)),\\[1ex]
\Lambda_{\eps} ((u_{\eps i}^\sigma)^{\pm}) \to u_i \mathcal{I} & \qtext{strongly % and a.e.
 in }L^2((0,T); L^{r}(\O)^{m\times m}).
\end{array}
\end{equation*}
In addition, for all $\eta\in L^{p'}((0,T);W^{1,p'}(\O))$ and $i=1,2$, we have
\begin{equation*}
\begin{array}{l}
\displaystyle\int _0^T\Big( \langle \pt u_{i}, \chi \rangle_{p'}
+ \big( u_{ i}  \big(a_{i1}\grad u_{ 1} +a_{i2}\grad
u_{2} +b_i q\big) + c_i\grad u_{ i} , \grad \chi \big)
\Big)
  = \\[1ex]
\hspace*{1cm} = 
\displaystyle\int _0^T\Big(  (\alpha_{i}u_{ i} - 
u_{ i} \big( \beta_{i1} \, u_{1} + \beta_{i2} \,
u_{2} \big) , \chi ) \Big) .
\end{array}
\end{equation*}
\end{lemma}

\bigskip

\subsection{Proof of Theorem~\ref{th:existence_particular}}

 We consider the perturbation $J_i^{BT,\delta}$ introduced in \fer{def.NLflowapp}, and recall that 
% $u^{(\delta)}=u_1^{(\delta)}+u_2^{(\delta)}$ satisfies Problem (P)$_\delta$, see \fer{prob:add.ec}-\fer{prob:add.ic}.
% and (H') holds, then 
% $u^{(\delta)}=u_1^{(\delta)}+u_2^{(\delta)}$ solves 
% the following problem in a weak sense
% \begin{align}
% & \pt u -\Div J^{(\delta)}(u)=f(u) && \qtext{in }Q_T=\O\times(0,T), 
%   \label{prob:add.ec}\\
% &  J^{(\delta)}(u)\cdot n =0 && \qtext{on }\Gamma_T=\partial\O\times(0,T),
% \label{prob:add.bc}\\
% & u(\cdot,0)=u_{0} && \qtext{in }\O,	\label{prob:add.ic} 
% \end{align}
% with $ u_0=u_{10}+u_{20}$, $J^{(\delta)}(u)=(a+\delta) u \grad u +b q u +c\grad u$ and 
% $f(u)= u (\alpha -\beta u)$. 
the following result is a consequence of Theorem~1 of \cite{chen04}.
\begin{lemma}\label{lem:existPeps}
\label{th:pert}
Assume (H'). 
Then there exists a weak solution $(u_1^{(\delta)},u_2^{(\delta)})$ of problem 
 (P)$_\delta$ in the following sense (for $i=1,2$):
\begin{itemize}
\item[(i)]  $u_{i}^{(\delta)}\geq 0$ satisfy the regularity properties
\begin{equation*}
\begin{array}{l}
u_{i}^{(\delta)}\in L^{\infty}(\Q)\cap L^{2}(0,T;H^{1}(\O  ))\cap
H^{1}(0,T;(H^{1}(\O  ))^{'}).
\end{array}
%\label{weak.regularity}
\end{equation*}
\item[(ii)]  For all $\vfi  \in L^{2} (0,T; H^{1}(\O ))$ we have 
\begin{eqnarray}
   \int_{0}^{T}\langle \pt u_{i}^{(\delta)} ,\vfi \rangle
    + \int_{\Q }  J_i^{(\delta)}(u_{1}^{(\delta)},u_{2}^{(\delta)}) \cdot\grad\vfi  
=\int_{\Q }f_i (u_{1}^{(\delta)},u_{2}^{(\delta)})~\vfi,  
\label{weak.pert} 
\end{eqnarray}
where
$\langle\cdot ,\cdot\rangle$ denotes the duality product of
$(H^{1}(\O ))^{'}\times H^{1 }(\O ).$
\item[(iii)]  The initial conditions \fer{eq:id} are satisfied
in the sense
\begin{equation*}
\begin{array}{l}
\lim%_{t \to 0}
\nor{ u_{i}^{(\delta)}(\cdot ,t)-u_{i0}}_{(H^{1}(\O ))^{'}}=0 
\end{array}
\quad \text{as}\quad t\rightarrow 0. 
%\label{weak.id}
\end{equation*}
\end{itemize}
In addition,
\begin{equation}
 \nor{u_{1}^{(\delta)}}_{L^\infty(\Q)} +  \nor{u_{2}^{(\delta)}}_{L^\infty(\Q)} + 
\nor{u_{1}^{(\delta)} +u_{2}^{(\delta)}}_{L^2(0,T;H^{2}(\O))}\leq C ,
\label{est:infh2}
\end{equation}
with $C$ independent of $\delta$.
\end{lemma}

\no\emph{Proof. }The proof is given in two steps.
Step one consists on showing the existence of 
solutions of problem (P)$_\delta$ using the ideas of the problem 
solved by Chen and J\"ungel \cite{chen04}, which strongly resembles ours. The only difference between both 
problems is in the definition of the diffusion matrices which, for the problem treated in \cite{chen04} 
is of the form
\[
A_1=\left(
 \begin{array}{cc}
  c_1+2a_{11}u_1+a_{12}u_2 & a_{12}u_1 \\
  a_{21}u_2  & c_2+2a_{22}u_2+a_{21}u_1
 \end{array}
\right),
\]
whereas for problem (P)$_\delta$ is given by
\begin{equation}
A_2=\left(
 \begin{array}{cc}
  c+(a+\delta)u_1+\frac{\delta}{2}u_2 & (a+\frac{\delta}{2})u_1 \\
  (a+\frac{\delta}{2})u_2  & c+(a+\delta)u_2+\frac{\delta}{2}u_1
 \end{array}
\right),
\label{A2:reg}
\end{equation}
which can not be recast in the form of $A_1$. However, as may be easily seen  in \cite{chen04}, 
this difference does not affect the proof as long as the matrix resulting from the change of unknowns 
$u_i=\exp(w_i)$ is symmetric and positive definite. And this is certainly the case since, rewriting 
the diffusion matrix obtained after this change of unknowns we get
\[
\tilde A_2=\left(
 \begin{array}{cc}
  c \mathrm{e}^{w_1}+(a+\delta)\mathrm{e}^{2w_1}+\frac{\delta}{2}\mathrm{e}^{w_1+w_2} & (a+\frac{\delta}{2})\mathrm{e}^{w_1+w_2} \\
  (a+\frac{\delta}{2})\mathrm{e}^{w_1+w_2}  & c\mathrm{e}^{w_2}+(a+\delta)\mathrm{e}^{2w_2}+\frac{\delta}{2}\mathrm{e}^{w_1+w_2}
 \end{array}
\right)
\]
which is  positive definite for all $\delta>0$. We may therefore adapt the results in \cite{chen04}
to obtain the existence of a solution of problem (P)$_\delta$, but in a weaker
sense than the notion of solution stated in Lemma \ref{th:pert}. 

The second step of the proof is intended to justify this point, and to this end we use that 
$u^{(\delta)}=u_{1}^{(\delta)}+u_{2}^{(\delta)}$ satisfies, in a weak sense, 
problem \fer{prob:add.ec}-\fer{prob:add.ic}.
Being this the case and recalling that assumptions (H') imply the uniform parabolicity of 
problem \fer{prob:add.ec}-\fer{prob:add.ic}, we may 
 apply  Theorem~3.1, Chapter V of \cite{LSU} to deduce uniform bounds for  
 $\nor{u^{(\delta)}}_{L^\infty(Q_T)}$  and $\nor{u^{(\delta)}}_{L^2(0,T;H^{2}(\O))}$. In particular, the 
 $L^\infty (\Q)$
bound on $u^{(\delta)}$ together with the non-negativity of $u_i^{(\delta)}$ obtained in \cite{chen04}
imply the uniform $L^\infty (\Q)$ bounds on $u_i^{(\delta)}$. In consequence, all the terms in the weak
formulation \fer{weak.pert} make sense for test functions in $L^2(0,T;H^1(\O))$. $\Box$

\bigskip

The proof of Theorem~\ref{th:existence_particular} is completed by 
passing to the limit $\delta\to0$.
Let $(u_{1}^{(\delta)},u_{2}^{(\delta)})$ be the solution of problem (P)$_\delta$
found in Lemma~\ref{th:pert}. 
As shown in \cite{chen04}, an entropy type inequality implies
\begin{equation}
\label{est:lin}
\sum_{i=1}^2 \int_{\Q} \big(2c \abs{\grad\sqrt{u_i^{(\delta)}}}^2 +\delta\abs{\grad u_i^{(\delta)}}^2 \big) \leq C, 
\end{equation}
with $C$ independent of $\delta$. In particular, the $L^\infty(Q_T)$ bound for $u_i^{(\delta)}$ found in Lemma~\ref{th:pert}
implies 
\begin{equation}
\label{bound.c}
 \int_{\Q} \abs{\grad u_i^{(\delta)}}^2  \leq 4 \nor{u_i^{(\delta)}}_{L^\infty(Q_T)} \int_{\Q} \abs{\grad\sqrt{u_i^{(\delta)}}}^2 \leq C,
\end{equation}
where $C$ capture several constants independent of $\delta$. However, observe that we only assume $c\geq 0$, so
bound \fer{bound.c} is irrelevant for the most interesting case of $c=0$.
From \fer{weak.pert} we deduce the following estimate for all $\vfi\in L^2(0,T;H^1(\O))$:
\begin{eqnarray*}
 && \int_{0}^{T}\langle \pt u_{i}^{(\delta)} ,\vfi \rangle \leq \nor{J_i^{(\delta)}(u_1^{(\delta)},u_2^{(\delta)})}_{L^2(Q_T)}  
 \nor{\grad\vfi}_{L^2(Q_T)} \\
&& \hspace{3cm} + \nor{f_i(u_1^{(\delta)},u_2^{(\delta)})}_{L^2(Q_T)}\nor{\vfi}_{L^2(Q_T)}.
\end{eqnarray*}
We have
\begin{eqnarray*}
&&\nor{J_i^{(\delta)}(u_1^{(\delta)},u_2^{(\delta)})}_{L^2(Q_T)}  \leq 
  a \nor{u_i^{(\delta)}}_{L^\infty(Q_T)} \nor{\grad (u_1^{(\delta)}+u_2^{(\delta)})}_{L^2(Q_T)} 
 \\
 &&  + \abs{b}  \nor{u_i^{(\delta)}}_{L^\infty(Q_T)} \nor{q}_{L^2(Q_T)} + c \nor{\grad u_i^{(\delta)}}_{L^2(Q_T)}
  +\delta \nor{u_i^{(\delta)}}_{L^\infty(Q_T)} \nor{\grad u_i^{(\delta)}}_{L^2(Q_T)} \\
 &&  +\frac{\delta}{2} \Big( \nor{u_1^{(\delta)}}_{L^\infty(Q_T)} \nor{\grad u_2^{(\delta)}}_{L^2(Q_T)}+ 
  \nor{u_2^{(\delta)}}_{L^\infty(Q_T)} \nor{\grad u_1^{(\delta)}}_{L^2(Q_T)} \Big).
\end{eqnarray*}
Using the $L^\infty(\Q)$ uniform estimates in \fer{est:infh2}, bound \fer{bound.c} 
 and assumptions  (H') we get
\begin{eqnarray*}
 \int_{0}^{T}\langle \pt u_{i}^{(\delta)} ,\vfi \rangle  \leq && 
 C_1\Big (1+ \nor{\grad (u_1^{(\delta)}+u_2^{(\delta)})}_{L^2(Q_T)} \Big) \nor{\grad\vfi}_{L^2(Q_T)} 
+ C_2\nor{\vfi}_{L^2(Q_T)} \\
&& +\delta C_3 \big( \nor{\grad u_1^{(\delta)}}_{L^2(Q_T)} +\nor{\grad u_2^{(\delta)}}_{L^2(Q_T)}\big)
 \nor{\grad\vfi}_{L^2(Q_T)},
\end{eqnarray*}
and from the $L^2(\Q)$ uniform estimates for $\grad (u_1^{(\delta)}+u_2^{(\delta)})$ 
in \fer{est:infh2} and estimate \fer{est:lin} we deduce
\begin{equation}
\nor{\pt u_i^{(\delta)}}_{L^2(0,T;(H^1(\O))'} \leq  C(1+\sqrt{\delta}).
\label{est:time}
\end{equation}
Thus, using \fer{est:infh2}, \fer{est:lin} and \fer{est:time}, we deduce the existence of subsequences
(not relabeled) and functions $u_i\in H^1(0,T;L^2(\O))\cap L^\infty(\Q)$ and 
$u\in H^1(0,T;L^2(\O))\cap L^2(0,T;H^2(\O))\cap L^\infty(\Q)$  such that
(see \cite{simon})
\begin{eqnarray}
&& \pt u_i^{(\delta)} \wto \pt u_i \qtext{weakly in }L^2(0,T;(H^1(\O))'),\nonumber\\
&& u_i^{(\delta)} \overset{*}{\wto} u_i \qtext{weakly * in }L^\infty(\Q), \label{debil}\\
&& \grad(u_{1}^{(\delta)}+u_{2}^{(\delta)}) \to \grad u \qtext{strongly in }L^2(\Q),\label{fuerte2}\\
&& u_{1}^{(\delta)}+u_{2}^{(\delta)} \to u \qtext{strongly in }L^2(\Q). \label{fuerte}
\end{eqnarray}
As a first observation, we may identify $u$ as $u_1+u_2$ due to \fer{debil} and \fer{fuerte}.
Using estimate \fer{est:lin} and the uniform $L^\infty(\Q)$ estimate of 
$u_i^{(\delta)}$ in \fer{est:infh2} we also deduce, for $i,j=1,2$ ,
\begin{eqnarray*}
 \delta \int_{\Q} u_i^{(\delta)}\grad u_j^{(\delta)} \cdot \grad \vfi  \leq 
\delta \nor{u_i^{(\delta)}}_{L^\infty(\Q))} \nor{\grad u_j^{(\delta)} }_{L^2(\Q))} \nor{\grad \vfi}_{L^2(\Q))} \leq 
C \sqrt{\delta} . 
\end{eqnarray*}
Finally, in the  passing to the limit $\delta\to 0$ in \fer{weak.pert} there are only two non-standard terms,
\[
 \int_{\Q}u_i^{(\delta)}\grad{(u_1^{(\delta)}+u_2^{(\delta)})}\cdot \grad\vfi 
\]
and 
\[
 \int_{\Q}f_i(u_1^{(\delta)},u_2^{(\delta)})\vfi = 
 \int_{\Q} u_i^{(\delta)}(\alpha -\beta (u_1^{(\delta)}+u_2^{(\delta)})\vfi  ,
\]
which converge to their corresponding limits in view of  \fer{debil}-\fer{fuerte}. $\Box$

\section{Conclusion}

We have shown that a natural election for cross-diffusion modeling, from the point of view of 
limit densities corresponding to systems of particles, is that introduced by Busenberg and
Travis \cite{busenberg83} from macroscopic ad-hoc considerations, in the discipline of population dynamics.
Although a rigorous deduction for boundary value problems has not been accomplished yet, 
the results for the Cauchy problem seems to point to the model considered in this article.
Mathematically, the problem of existence of solutions has two cases. The first is the case in
which the system matrix $(a_{ij})$ is positive definite, for which we have given a rather general proof 
based on previous results for the Shigesada et al. model \cite{shigesada79}. The second, is the
case in which this matrix is only positive semi-definite. We have given a partial result of 
existence of solutions which generalizes previous results based on the solution construction 
by Lagrangian flows. In this case, the problem is specially interesting for segregated initial data,
giving rise to the contact inhibition problem arising from tumor modeling. After checking the qualitative
similarities, from a numerical simulation point of view, between the BT and the SKT models when 
the problem is parabolic (positive definite matrix), we have reviewed several situations in which the presumably 
non-parabolic problem (positive semi-definite matrix) gives
rise to discontinuous solutions. We have also performed simulations out of the range of the assumptions
for the existence proof, showing that they seem to be just technical restrictions. In future work, 
we shall investigate the possibilities of broadening such conditions.

\section*{Acknowledgements}

The authors thank to the anonymous reviewers for their comments and observations. They have
contributed to the improvement of our work.

\end{document}